\documentclass[a4paper, 12pt]{article}
\usepackage[a4paper, left=2cm, right=1cm]{geometry}
\usepackage{multirow}
\usepackage{array}
\usepackage{color}
\usepackage{amsmath,amssymb,latexsym}
\usepackage[mathscr]{eucal}
\usepackage{arydshln}
\newcommand{\R}{\mathbb{R}}
\newcommand{\Z}{\mathbb{Z}}

\newcommand{\N}{\mathbb{N}}

\newcommand{\C}{\mathbb{C}}
\newcommand{\beq}{\begin{equation}}
\newcommand{\eeq}{\end{equation}}
\newcommand{\bea}{\begin{eqnarray}}
\newcommand{\eea}{\end{eqnarray}}
\newcommand{\bean}{\begin{eqnarray*}}
\newcommand{\eean}{\end{eqnarray*}}

\newtheorem{theorem}{Theorem}
\newtheorem{lemma}[theorem]{Lemma}
\newtheorem{coro}[theorem]{Corollary}
\newtheorem{rem}[theorem]{Remark}
\newtheorem{pro}[theorem]{Proposition}

\newenvironment{proof}%
{\par\noindent\emph{Proof:\ }}%
{\ \hfill ~\rule{2mm}{2mm}\par\bigskip}
{\par\noindent\textbf{Remark:\ }}%
{\ \hfill \par\bigskip}

\def\txx{\tilde{x}}
\def\tyy{\tilde{y}}
\def\tzz{\tilde{z}}
\def\ta{\tilde{a}}
\def\tgam{\tilde{\gamma}}
 \def\ta{\tilde{a}}
\def\dd{\tilde{d}}

\def\hxx{\hat{x}}
\def\hyy{\hat{y}}
\def\hzz{\hat{z}}
\def\hgam{\hat{\gamma}}
\def\ha{\hat{a}}
\def\hc{\hat{c}}
\def\heps{\hat{\epsilon}}
\def\ee{\hat{d}}

\def\bxx{\bar{x}}
\def\byy{\bar{y}}
\def\bzz{\bar{z}}

\def\bgam{\bar{\gamma}}

\def\mut{{{\mu}}}
\def\muu{{{\nu}}}   
\def\mud{{{\nu}}}   
\def\muh{{{\hat{\nu}}}}   
\def\hd{{{h}}}   
\def\mutd{{{\tilde{\nu}}}}   
\def\htd{{{\tilde{h}}}}   
\def\hmud{{{\hat{\mu}}}}
 \def\landa{{{\lambda}}}
\def\tlambda{{\tilde\landa}}
\def\hlambda{\hat{\landa}}
 \def\CCp{\hat{\mathcal{C}}}
 \def\DDp{\hat{\mathcal{D}}}
 \def\FFp{\hat{\mathcal{F}}}
 \def\CCt{{\mathcal{C}}}
 \def\DDt{\mathcal{D}}
 \def\FFt{\mathcal{F}}
 \def\GGt{\mathcal{G}}
 \def\g{{{h}} } 
 \def\tg{{\tilde{h}}}   
 \def\tmu{{{\tilde\mu}}}
 \def\varsig{{{f}}}
\def\hpp{\mathfrak{h}}
\def\fpp{\mathfrak{f}}
\def\fantasm{\phantom{\Bigg|}}

\def\Gh{\widehat{\mathbf{G}}}

\makeatletter
\newcolumntype{"}{@{\hskip\tabcolsep\vrule width 2pt\hskip\tabcolsep}}
\def\hlinewd#1{%
  \noalign{\ifnum0=`}\fi\hrule \@height #1 \futurelet
   \reserved@a\@xhline}
\makeatother
\def\t{{\mathbf t}}
\def\ttr{{\mathbf t}}
\def\tpl{{\hat{\ttr}}}
\def\nablatr{{\nabla}}
\def\nablapl{{\hat{\nabla}}}
\def\Opl{{\ell}}
 \def\Oplt{{\ell}}
 \def\Opltc{{\ell}^{c}}
 \def\Oplp{\hat{\ell}}
 \def\qh{\mathscr{P}}
 \def\QH{\mathcal{Q}}
 \def\spacE{\hat{\mathcal E}}
 \def\spacEt{{\mathcal E}}
 
 \def\spacKt{{\mathcal K}}
 \def\Deltat{{\Delta}}
 \def\Deltap{\hat{\Delta}}
 \def\spacV{{\mathcal V}}
 \def\spactV{\widetilde{\mathcal V}}

 \def\Lconj{\mathbf{L}}
 \def\Lequiv{{\mathcal{L}}}

 \def\x{{\mathbf x}}
 \def\y{{\mathbf y}}
 \def\F{{\mathbf F}}
 \def\G{{\mathbf G}}
 \def\P{{\mathbf P}}
 \def\X{{\mathbf X}}
 \def\r{{\mathtt r}}
 \def\c{{\mathtt c}}
 \def\zero{{\mathbf 0}}
\def\tx{\widetilde{x}}
\def\ty{\widetilde{y}}
\def\tz{\widetilde{z}}
\def\tP{{\widetilde{\mathbf P}}}
\def\tF{{\widetilde{\mathbf F}}}
\def\FZERO{{\mathbf{D}_0}}
\def\spande{{\mathrm{Span}}}

\def\diverg{{\mathrm{div}}}
\def\Cor{{\mathrm{Cor}}}
\def\Ker{{\mathrm{Ker}}}
\def\Range{{\mathrm{Range}}}
\def\Proy{{\mathrm{Proj}}}
\def\bracket#1{\left[ #1 \right]}
\def\parent#1{\left( #1 \right)}
\def\floor#1{\left\lfloor {#1}\right\rfloor}
\def\llave#1{\left\{ #1 \right\}}

\def\vedo#1#2{ \begin{pmatrix} #1 \\  #2 \end{pmatrix}  }
\def\vetre#1#2#3{ \begin{pmatrix} #1 \\ #2 \\  #3 \end{pmatrix}  }
\def\vedot#1#2{\left(  #1 , #2 \right)^{T}}
\def\vetret#1#2#3{ \left(  #1 , #2 , #3  \right)^T }
\def\vedor#1#2{ \begin{pmatrix} #1 \\  \hdashline[2pt/2pt] #2 \end{pmatrix}  }
\def\vetrer#1#2#3{ \begin{pmatrix} #1 \\ #2 \\ \hdashline[2pt/2pt] #3 \end{pmatrix}  }

\def\fracp#1#2{{\textstyle{\frac{#1}{#2}}}}
\date{\today}
\title{{
Orbital Normal Forms for a Class of  Three-Dimensional Systems with an Application to Hopf-Zero Bifurcation Analysis of Fitzhugh--Nagumo System}}

\author{
A. Algaba$^\dagger$, N. Fuentes$^\dagger$, E. Gamero$^\ddagger$, C. Garc\'{\i}a$^\dagger$
\\
$^\dagger$Dept. Integrated Sciences. Faculty of Experimental Sciences\\ University of Huelva,
Spain
\\
$^\ddagger$ Dept. Applied  Mathematics II, E. T. S. I.\\  University of
Sevilla, Spain}
\begin{document}
\bibliographystyle{plain}

\maketitle

\abstract{{
We consider a class of three-dimensional systems having an equilibrium point at the origin, whose  principal part is of the form 
$\vetret{-\frac{\partial {\hpp}}{\partial y}(x,y)}
              { \frac{\partial {\hpp}}{\partial x}(x,y)}  
             { \fpp(x,y) }$.
This principal part, which has zero divergence and does not depend on the third variable $z$,  is the coupling of a   planar Hamiltonian vector field $\X_{\hpp(x,y)}:=\vedot{-\frac{\partial {\hpp}}{\partial y}(x,y)}
              { \frac{\partial {\hpp}}{\partial x}(x,y)} $ 
with a one-dimensional system.

We analyze the quasi-homogeneous orbital normal forms for this kind of systems, by  introducing a new splitting for quasi-homogeneous three-dimensional vector fields.
The obtained results  are applied to the nondegenerate Hopf-zero singularity that falls into this kind of systems.
Beyond the Hopf-zero normal form, a parametric normal form is obtained, and the analytic expressions for the normal form coefficients are provided.
Finally, the results  are applied to a case of the three-dimensional Fitzhugh-Nagumo system.
}

\section{Introduction} 

The underlying idea 
of the classical Poincar\'e normal form theory  for equilibria of smooth autonomous vector fields
is to use  near-identity transformations to remove degree by degree the nonessential nonlinear terms
in the Taylor expansion of the vector field, so that  in the normalized vector field only  resonant nonlinear terms remain.
These resonant terms are completely determined by the linear part.
Generically, this normal form is not unique because the homological equation (which determines the simplifications that can be achieved) can have infinitely many solutions,
and it is possible to use the elements of the kernel of the homological operator for obtaining simplified normal forms.

{
There has been numerous contributions in the study of simplest normal forms, mainly in the case of planar systems. 
The normal form classification for the Takens-Bogdanov singularity has been completely solved in \cite{StroZol15}. For the case of the Hopf singularity see, e.g. \cite{Algaba99,YuLeung03}. A  case of planar system with null linear part has been considered in \cite{Stroz17}. The normal form problem becomes significantly more complicated for three-dimensional systems.
Most of the contributions in this case corresponds to the Hopf-zero singularity, see, e.g., \cite{YuYuan00}, and also to  nilpotent three-dimensional singularities (see \cite{AlgCL,Li17,Li14}.

In the our normal form analysis}, we select a type $\t$ and use here QH (quasi-homogeneous) expansions  stead of Taylor expansions.
As above, the idea is to simplify such expression degree by degree in the QH ordering by using 
near-identity transformations to annihilate nonessential terms in the local dynamical behavior of the system.
The simplifications are determined by the \emph{principal part} of the vector field with respect to the type $\t$ (the lowest-degree  term in the QH expansion of the vector field).
Changing the type $\t$  leads to  different  principal parts and therefore, to different normal forms.
Then, it is very important to select adequately the type in order to obtain normal forms as simple as possible.
Ideally, if the type is selected such that the homological equation has a unique solution, then we reach the unique normal form (i.e., the simplest one).

For instance, the classical normal form for the nondegenerate Takens-Bogdanov singularity
is based on the linear part $\vedot y0$, but it is more convenient to use a QH expansion
with principal part $\vedot y{x^2}$ because one obtain simpler normal forms (see \cite{Algaba02}).

In the same way, in the analysis of three-dimensional systems undergoing a nondegenerate Hopf-zero singularity,
it is usual to use Taylor expansions to determine normal forms based on the linear part $\vetret{-y}x0$.
Nevertheless, we will see in this paper that it is more convenient to use QH expansions having principal part
$\vetret{-y}x{x^2+y^2}$ in order to get simplified normal forms.

{
We notice that the principal part  $\vetret{-y}x{x^2+y^2}$ for the  nondegenerate Hopf-zero singularity has the structure that we consider in this paper. There are other interesting cases that also fall into the class we are considering; for instance the non-semisimple nilpotent singulatiry in dimension  three corresponding to the coupling of Takens-Bogdanov and saddle-node singularities,  whose principal part is of the form $\vetret y{x^2}{ax^3+by^2}$,   or different cases of three-dimensional systems with null linear part.

In these cases, the use of adequate QH expansions provides simplified normal forms. In fact, it is the most natural ordering because the quasi-homogeneity type determines the rescaling if we use blow-up techniques.

}

A major feature that must be kept in mind when obtaining normal forms is related to its applicability.
In fact, it is not so much a question of obtaining \emph{simplified} normal forms as of obtaining \emph{useful} normal forms.
In this regard,  we obtain in this paper an orbital normal form for the nondegenerate Hopf-zero singularity that  is not the simplest one (a hypernormal form procedure should be carried out in order to get further simplifications),
but it is appropriate to analyze  the integrability problem for such singularity, as we  show in \cite{AlgHZInteg}.

A basic idea in the study of normal forms is the analysis of the Lie product.
In the case of planar systems, the use of a conservative-dissipative splitting  has proved to be very effective for analyzing not only normal forms, 
but also for some classical problems in this context. 
Namely, the quoted splitting has been used in the study of the analytic integrability problem 
(see \cite{AlgabaNonlinearity09,AlGarGi,AGGintnilp16}),  in the analysis of the center problem (see \cite{AlGarGa}) 
and also to determine the existence of inverse integrating factors (see \cite{Algaba12,Algaba122}) 
or  symmetries  (see  \cite{AGGC}).

Here, we  present a new  splitting for three-dimensional
QH vector fields that generalizes the planar splitting just mentioned.
This splitting is very useful for the analysis of the normal forms and  it can be used
in the study of the analytical integrability problem in some cases of three-dimensional systems
(see \cite{AlgHZInteg}).

In this paper, we deal with a large class of three-dimensional systems of the form
\bea
\nonumber
\dot x &=& F(x,y,z), \\
\label{1}
\dot y &=& G(x,y,z), \\
\nonumber
\dot z &=& H(x,y,z), 
\eea
{
having an equilibrium at the origin, i.e., 
$$F(0,0,0)=0, \ G(0,0,0)=0, \ H(0,0,0)=0.$$
We assume that the principal part  of system (\ref{1}) with respect to some type $\t$ is given by:
\beq
\label{2}
\vetre{-\frac{\partial {\hpp}}{\partial y}(x,y)}
              { \frac{\partial {\hpp}}{\partial x}(x,y)}  
             { \fpp(x,y) },
\eeq
where $\hpp(0,0)=\fpp(0,0)=0$.
Principal part here refers to the lowest degree terms appearing in the expansion of the system around the singular point, i.e.,
\bean
F(x,y,z)&=& {-\frac{\partial {\hpp}}{\partial y}(x,y)} + \mathrm{h.o.t.}, \\
G(x,y,z)&=&  { \frac{\partial {\hpp}}{\partial x}(x,y)} + \mathrm{h.o.t.},  \\
H(x,y,z)&=&  { \fpp(x,y) }+ \mathrm{h.o.t.}, 
\eean
where h.o.t. denotes higher QH order terms.
}

We notice that  principal part (\ref{2}) has zero divergence and only depends on two variables (it is independent on $z$).
Also, it is the coupling of a   planar Hamiltonian vector field associated with the Hamiltonian $\hpp(x,y)$ and  a one-dimensional system.

The above structure appears in a wide class of systems which include, as particular cases,
the normal form of nondegenerate Hopf-zero singularity commented before. 
Also, the degeneracy corresponding to  the interaction of  Takens-Bogdanov and  fold bifurcations  falls in this case.
{
In \cite{AlgCL}, the normal form for this singularity, that is associated  with a triple zero eigenvalue with geometric multiplicity two, has been obtained.
In fact, in  \cite{AlgCL} some results that are obtained here are presented without proof. 
Moreover, here we obtain normal forms under conjugation as well as under orbital equivalence.
Instead, in   \cite{AlgCL}  only normal forms under conjugation are considered.}

This paper is organized as follows.
In the next section, we present a brief review of the essential concepts about normal forms
and its characterization  using QH expansions.
Later, in Section \ref{TechnicalResults} we provide a new splitting for QH
three-dimensional vector fields which is useful for the characterization of  normal forms
for system (\ref{1}) having principal part (\ref{2}). 
{
We carry out this in Section \ref{NormalForm}; see  Theorems \ref{NormaFormFinal-}  and \ref{tecorlylc}.}
Finally,  in Section  \ref{Hopzero} the results are applied to the case of the nondegenerate Hopf-zero singularity. 
Also, we obtain a parametric normal form, and we present the analytic expressions for the normal form coefficients.
As a case study, we consider the three-dimensional Fitzhugh-Nagumo system with some parameters values where a Hopf-zero degeneracy takes place. We use the parametric normal form results  to obtain information about periodic behavior for the {
Fitzhugh-Nagumo}  system.

\section{Basic concepts on Normal Form Theory}

As shown in \cite{Algaba02,Algaba04}, the use of  QH expansions of  vector fields in the  study of normal forms
has important advantages because in the QH ordering linear and nonlinear terms are mixed.
If one wants to use the classical Taylor expansions, it is enough to select the unity type.
We include below a short summary of the results about  QH normal forms, that  can be found in more detail in {
\cite{Algaba02,Algaba04}.}

A type $\t=(t_1, \dots, t_n)\in\N^n$ is a $n$-index that determines the weights of the variables
(for this reason,  QH polynomials are also called   weighted-homogeneous polynomials).
We use standard multi-index notation; for instance $|\t|:=t_1+\cdots+t_n$ denotes the modulus of $\t$.

A scalar function $f$ of $n$ variables is QH of type
$\t\in\N^n$ and degree $k$ if
$f(\epsilon^{t_1} x_1, \dots , \epsilon^{t_n} x_n )= \epsilon^k f(x_1, \dots, x_n)$.
The vector space of QH functions of type $\t$ and degree $k$ is denoted by $\qh^{\t}_k$.

A vector field $\F=(F_1, \dots, F_n)^T$ is   QH of type $\t$ and degree $k$ if
$F_j\in \qh^{\t}_{k+t_j}$ for each $j$.
We  denote $\QH_k^{\t}$ the vector space of QH vector fields of type $\t$ and degree $k$.

Let us consider an analytic  autonomous system $\dot{\x}=\F(\x)$, with $\x=(x_{1},\dots,x_{n})\in\R^{n}$
having an equilibrium point at the origin (i.e., $\F(\zero)=\zero$).
The  vector field $\F$ can always be  written as the sum of QH terms of type $\t$:
 \beq \label{sisprinc}
\dot{\x} = \F(\x)  = \F_{r}(\x) + \F_{r+1}(\x) + \cdots ,
 \eeq
where $\F_{k}\in\QH^{\t}_{k}$ for all $k$.
The lowest-degree QH $\F_{r}\not\equiv \zero$ (being $r \in \Z$) is  called the
\emph{principal part} of $\F$ with respect to the type $\t$.

The underlying idea in the normal form theory is to remove unessential terms in the analytical expression of a vector field.
The simplifications  can be obtained through transformations in the state variables
(\emph{normal forms under conjugation}). 
They are determined by the homological operator, that depends on the principal part $\F_r$ and  is defined 
by means of the  Lie product of two vector fields: {
$\bracket{\P,\F}:= D\P \, \F - D\F \, \P$, where $D$ stands for the differential with respect to the state variables.}

The simplification procedure can also use time-reparametrizations (\emph{normal forms under orbital equivalence}
or \emph{orbital normal forms}).
In this case, if we reparametrize the time by $\frac{dt}{dT}=\mut(\x)$, with $\mut(\zero)=1$, then the transformed vector field
agrees with the original one multiplied by $\mut$.
To combine the effect of time-reparametrizations and transformations in the state
variables, it is enough to use that {
$ \bracket{\P,\mut\F }=\mut\, \bracket{\P, \F}-(\nabla \mut\cdot \P)\F$,
for any smooth scalar function $\mut$ and vector fields  $\P$, $\F$ (see \cite{GazorYu12}).
Here, we assume that the time rescaling is first used and then the changes of state variables are applied.}

Next, we enumerate  several properties related to the use of QH vector fields
(their proofs can be found, e.g., in \cite{Algaba02}). 
\begin{itemize}
\item
If $\P \in \QH_{k}^{\t}$, $\F \in \QH_{l}^{\t}$, then $[\P,\F]  \in \QH_{k+l}^{\t}$.
\item
If $\mut\in\qh_{k}^{\t} $ and $ \F\in \QH_l^{\t}$, then $\mut\,\F\in \QH_{k+l}^{\t}$.
\item
If $\mut\in\qh_{k}^{\t}$ and $\F\in \QH_l^{\t}$, then $\nabla \mut\cdot \F\in \qh_{k+l}^{\t}$.
\item
The  Hamiltonian vector field associated with the Hamiltonian $\hpp(x,y)$ is  
$$\X_{\hpp} = \vedo{-\frac{\partial \hpp}{\partial y} }{ \frac{\partial \hpp}{\partial x} }.$$
We observe that $\X_\hpp\in\QH_k^{\t}$ if, and only if, $\hpp\in\qh_{k+|\t|}^{\t}$.
 \item
 If $\F \in \QH_k^\t$, then its divergence  satisfies $\diverg\parent{ \F } \in\qh_k^\t$.
 \item
 The  wedge product (see \cite{Guckenheimer83b}) of two planar  vector fields
 $\F=(P,Q)^{T} \in \QH_{k}^\t$, $\G=(R,S)^{T} \in \QH_{l}^\t$
 satisfies  $ \F \wedge \G =  P S - Q R  \in \qh_{k+l+|\t|}^{\t}$.
 \end{itemize}

\subsection{Normal forms under conjugation}
\label{SQHNFC}

The basic idea in the classical study of normal forms for system (\ref{sisprinc}) is the analysis
of the effect of a near-identity transformation $\x=\y+\P_k(\y)$, where $\P_k\in\QH_{k}^{\t}$
with $k\geq 1$.
The transformed system is
 $$
\dot{\y} = \G(\y) = \parent{ I + D \P_k (\y) }^{-1} \sum_{j\geq 0} \F_{r+j}\parent{ \y + \P_k (\y) }.
 $$
It can be shown that the QH terms of this transformed system agree with the
original ones up to degree $r+k-1$, and at degree $r+k$ we have:
$$
\G_{r+k} = \F_{r+k} - \parent{ D \P_k\, \F_{r} - D \F_{r}\, \P_k } =
\F_{r+k}- \bracket{\P_k, \F_{r}}=\F_{r+k} - \Lconj_{r+k}(\P_k),
$$
where we have introduced the \emph{homological operator}:
 \bean
\Lconj_{r+k} &:& \QH_k^{\t} \longrightarrow \QH_{r+k}^{\t}
\\ &&
\P_k \rightarrow \Lconj_{r+k} (\P_k)=\bracket{\P_k, \F_{r}}.
 \eean
Hence, for each $k\geq 1$ we can annihilate  in $\G_{r+k}$ the part belonging to the range of the
linear operator $\Lconj_{r+k}$, by selecting $\P_k$ adequately.
Namely, it is enough to write $\F_{r+k}=\F^\r_{r+k}+\F^\c_{r+k}$ where
$\F^\r_{r+k}\in\Range\parent{\Lconj_{r+k}}$ and $\F^\c_{r+k}\in\Cor\parent{\Lconj_{r+k}}$: a complementary subspace to the range of $\Lconj_{r+k}$.
By selecting $\P_k$ satisfying the homological equation $\Lconj_{r+k}(\P_k)=\F^\r_{r+k}$,
we achieve $\G_{r+k}=\F_{r+k} - \F^\r_{r+k} = \F^\c_{r+k}$, and we say that this term has been reduced
to normal form.
Performing  this procedure for $k=1,2,\dots$, we obtain the following result:
\begin{theorem} \label{NF}
System (\ref{sisprinc}) can be formally reduced to normal form:
$$
\dot{\y} = \G(\y) = \G_r (\y) + \G_{r+1}\parent{ \y  }+\cdots,
 $$
where $\G_r=\F_r$ and $\G_{r+k}\in\Cor\parent{\Lconj_{r+k}}$ for all $k\geq 1$,
by a sequence of near identity transformations.
\end{theorem}

\subsection{Normal forms under orbital equivalence}
 \label{SQHNFE}

The  simplifications of the  \emph{Normal Form Theorem} \ref{NF}  can be improved by including the effect of
time-reparametrizations in the simplification procedure (see \cite{Algaba04}).
In this case, for each $k\geq 1$, we first reparametrize the time by $\frac{dt}{dT}=1+\mut_k(\x)$, with
$\mut_k\in \qh^{\t}_k$ and later, we use  a near-identity transformation $\x=\y+\P_k(\y)$.
It is immediate to show that again the transformed system $\dot \y=\G(\y)$ agrees with the original one up to degree
$r+k-1$, and the $(r+k)$-degree term is:
 $$
\G_{r+k} = \F_{r+k} + \mut_k \, \F_{r} - \Lconj_{r+k}(\P_k) = \F_{r+k} - \Lequiv_{r+k}\parent{\P_k, \mut_k},
 $$
where we have introduced the \emph{orbital homological} operator:
 \bea
\Lequiv_{r+k}  &:& \QH_k^{\t} \times \qh_k^{\t} \longrightarrow \QH_{r+k}^{\t} \label{OpHomEquiv}
 \\ &&
\parent{\P_k, \mut_k} \rightarrow \Lequiv_{r+k}(\P_k, \mut_k) = \bracket{\P_k, \F_{r}}- \mut_k\, \F_{r} .\nonumber
 \eea
Reasoning as above, it is possible to choose  $(\P_k, \mut_k)$
adequately in order to simplify the $(r+k)$-degree QH term in system (\ref{sisprinc}),
by annihilating the part belonging to the range of   $\Lequiv_{r+k}$
(in this way, we  achieve that $\F_{r+k}$ belongs to a complement to the range of $\Lequiv_{r+k}$).
When this has been done, we say that the corresponding term has been reduced to orbital normal form.
By performing the procedure for $k=1,2,\dots$, we obtain the \emph{Orbital Normal Form Theorem}:
\begin{theorem} \label{NFE}
System (\ref{sisprinc}) can be formally reduced to orbital normal form:
$$
\dot{\y} = \G(\y) = \G_r (\y) + \G_{r+1}\parent{ \y  }+\cdots,
$$
where $\G_r=\F_r$ and $\G_{r+k}\in\Cor\parent{\Lequiv_{r+k}}$ for all $k\geq 1$,  by a sequence of
time-reparametrizations and near identity transformations.
\end{theorem}

{
We notice that the orbital normal form of Theorem \ref{NFE} might not be simplest one, and a further hypernormalization procedure should be performed, but this is beyond the scope of the present paper.}

The homological operators $\Lconj_{r+k}$, $\Lequiv_{r+k}$ depend on the principal part ${\F_{r}}$ which, in turn,  depends on
the type $\t$ chosen.
For this reason, the selection of the unit type $\t=(1,\dots,1)$ (that corresponds to the classical normal form based on Taylor expansions)  is not always the best:
a smart selection could provide simpler normal forms, as we will see later for  the Hopf-zero singularity considered in
Section \ref{Hopzero}.

{
The orbital normal form analysis is more difficult than the normal form computation using only changes of state variables.
This is the reason that there are considerably less research works on orbital normal forms.}
The analysis of the orbital homological operator  $\Lequiv_{r+k}$ is more involved that the
one of conjugation because some part of $\F_{r+k}$ can be annihilated with transformations 
in the state variables (through $\P_k$)  as well as with time-reparametrizations 
(through $\mut_k$). 
In \cite{Algaba02},  it is shown that the simplifications can be reached by taking $\mut_{k}\in \Cor\parent{\Opl_k}$, 
a complement to the range  of the Lie derivative operator associated with the principal part of the 
vector field (\ref{sisprinc}), which is defined by:
\bea
\Opl_{k} &:& \qh_{k-r}^\t \longrightarrow \qh_{k}^\t \label{Opele}
\\
&& \,  \muu_{k-r}\longrightarrow\Opl_{k}\parent{ \muu_{k-r}}=\nabla  \muu_{k-r}\cdot\F_{r}.
\nonumber
\eea

Hence, it is enough to consider the orbital homological operator with a restricted domain:
$$
\Lequiv_{r+k} : \QH_k^{\ttr} \times \Cor(\Oplt_{k}) \longrightarrow \QH_{r+k}^{\ttr}.
$$

\section{{
Some splittings for quasi-homogeneous vector fields}}
\label{TechnicalResults}

We are concerned with three-dimensional systems expanded in quasi-homogenous terms of some fixed type
$\t=(t_1, t_2, t_3)$:
\beq\label{sistprinspecifically}
\dot{\x}=\F(\x) := \F_r(x,y) + \F_{r+1}(\x) + \cdots,
\eeq
where  $\x=(x,y,z)\in\R^3$,  and the principal part, which  does not depend on $z$, is of  the form:
\beq\label{ppalpart}
\F_{r}(x,y) = \vetre{-\frac{\partial {\hpp}}{\partial y}(x,y)}{ \frac{\partial {\hpp}}{\partial x}(x,y)}{\fpp(x,y)}   =  \vedor{ \X_{\hpp}(x,y) }{ \fpp(x,y) },
\eeq
and the higher order QH terms $\F_{r+1},\dots$
depend on the three variables.

Along this paper, we use QH functions depending on  two and three variables.
To  distinguish them, we use a hat in the planar context.
For instance, given the type  $\ttr=(t_1, t_2, t_3)$, we  denote  $\tpl=(t_1, t_2)$.
 This type will appear when dealing with  functions depending on two variables:
notice that  the principal part $\F_{r}\in\QH_{r}^{\ttr}$, whereas
$$\X_{\hpp}\in\QH_{r}^{\tpl} \
\parent{\textrm{i.e., }\hpp\in \qh_{r+|\tpl|}^{\tpl}},
\textrm{  and } \fpp \in \qh_{r+t_3}^{\tpl}.
$$

In the same way, we  use the  three-dimensional Lie  derivative operator $\Oplt_{k}$ associated with $\F_{r}$
(see (\ref{Opele})),
 as well as the  Lie  derivative operator associated with the Hamiltonian planar vector field $\X_{\hpp}$, which we denoted as
\bea \label{derLie2-dim}
\Oplp_{k} &:& \qh_{k-r}^{\tpl} \longrightarrow \qh_{k}^{\tpl}
\\ && \nonumber
\, \muh_{k-r} \longrightarrow
\Oplp_{k}\parent{\muh_{k-r}}=
\nablapl \muh_{k-r}\cdot\X_{\hpp}, 
\eea
where  the planar gradient is denoted by
$\nablapl:=\parent{\frac{\partial}{\partial x},\frac{\partial}{\partial y}}^T$, 
whereas the three-dimensional gradient is
$\nablatr:=\parent{\frac{\partial}{\partial x},\frac{\partial}{\partial y},\frac{\partial}{\partial z}}^T$.

Next result shows that the principal part (\ref{ppalpart}) can be adequately simplified.
\begin{pro}\label{finCorl}
{
Let us consider system (\ref{sistprinspecifically}) with the principal part  given in (\ref{ppalpart}). Then, there exists a near-identity transformation such that    
$$
\fpp\in\Cor\parent{\Oplp_{r+t_3}}.
$$
Here, $\Cor\parent{\Oplp_{r+t_3}}$  a complement to the range of the Lie derivative operator $\Oplp_{r+t_3}$.
}
\end{pro}

\begin{proof}
Let us consider the system (\ref{sistprinspecifically}), and perform the change of variables
$\tx=x$, $\ty=y$, $\tz=z-\muu_{t_3}(x,y)$, with $\muu_{t_3}\in\qh_{t_3}^{\tpl}$.
It is easy to show that the transformed system has the following principal part:
$$\tF_{r} = \vedor{ \X_{\hpp} }{ \fpp - \nablapl \muu_{t_3}\cdot\X_{\hpp} } = \vedor{ \X_{\hpp} }{ \fpp - \Oplp_{r+t_3}\parent{ \muu_{t_3} } }.$$
To complete the proof, it is enough to select
 $\muu_{t_3}$ adequately in order to annihilate 
the projection of $\fpp$ onto the range of $\Oplp_{r+t_3}$.
\end{proof}

From now on,  we suppose that  the principal part $\F_r$ given in (\ref{ppalpart}) has been simplified according to the above proposition. 
Moreover, we assume that the factorization of $\hpp$ on $\C[x,y]$ only has simple factors.
This hypothesis warrants that,  we have
\beq\label{Corciclico}
\Cor\parent{\Oplp_{r+k+|\tpl|}}=\hpp\Cor\parent{\Oplp_{k}},
\eeq
(see  Proposition 3.18 of \cite{AlgabaNonlinearity09}). 
Equality (\ref{Corciclico})  allows to characterize the complements $\Cor\parent{\Oplp_{k}}$ (for all $k$)
in a finite number of steps.

\subsection{{
A splitting for quasi-homogeneous planar vector fields}}
\label{DecomposPla}

A  conservative-dissipative splitting of
QH planar vector fields has been used in \cite{AlgabaNonlinearity09} 
in the study of the  integrability problem.
{
Similar instances have been used in \cite{ChassIchi19} for the use of Melnikov integrals and in \cite{PalacYang00}
for developing efficient ODE solvers.
An alternative for the conservative-dissipative splitting of the normal forms is to use the Lie algebra characterization of singular vector fields, see \cite{GazShoPrepr}.

The conservative-dissipative splitting is introduced in the following proposition.}

\begin{pro}\label{con-dis}
Let us consider $\P_k\in\QH_{k}^{\tpl}$ and denote $\FZERO=\vedot{t_1 x}{t_2 y} \in \QH_{0}^{\tpl}$.
Then, there exist unique QH polynomials
$\hd_{k+|\tpl|}\in\qh_{k+|\tpl|}^{\tpl}$ and $\mud_k \in \qh_{k}^{\tpl}$  such that:
\beq\label{condis}
\P_k=\X_{\hd_{k+|\tpl|}} + \mud_k \, \FZERO.
\eeq
Moreover,
$\hd_{k+|\tpl|} = \frac{1}{k+|\tpl|}\FZERO\wedge\P_{k}$ and 
$\mud_k = \frac{1}{k+|\tpl|}\diverg(\P_{k})$.
\end{pro}
Several  results regarding the behavior
of the above splitting for QH functions are also derived in \cite{AlgabaNonlinearity09}.
In particular, the following result is proved.
\begin{lemma}\label{Oper}
\begin{itemize}
\item[(a)]
If   $h_k \in \qh_{k}^\tpl$ and $h^*_{l+|\tpl|} \in \qh_{l+|\tpl|}^\tpl$, then
$\bracket{ \X_{h_k},\X_{h^*_{l+|\tpl|}} } = \X_{\htd_{k+l}}$, where 
$$\htd_{k+l}=\nablapl h_k \cdot \X_{h^*_{l+|\tpl|}}\in \qh_{k+l}^\tpl.$$
\item[(b)]
Let us consider $\mut_k \in \qh_{k}^\tpl$ and $h^*_{l+|\tpl|} \in \qh_{l+|\tpl|}^\tpl$. Then,
$\mut_k\X_{h^*_{l+|\tpl|}} = \X_{\htd_{k+l+|\tpl|}}$ + ${\mutd_{k+l}} \, \FZERO$, with
\bean
\htd_{k+l+|\tpl|}&=&\frac{l+|\tpl|}{k+l+|\tpl|}\mut_k h^*_{l+|\tpl|}\in\qh_{k+l+|\tpl|}^{\tpl}, \textrm{ and }
\\
{\mutd_{k+l}}&=&\frac{1}{k+l+|\tpl|}\nablapl \mut_k\cdot \X_{h^*_{l+|\tpl|}} \in \qh_{k+l}^{\tpl}.
\eean
\end{itemize}
\end{lemma}

Let us  define the following subspace, which depends on the Hamiltonian function  $\hpp\in \qh_{r+|\tpl|}^{\tpl}$
given in the principal part (\ref{ppalpart}):
\beq\label{Delta}
\hpp\qh_{k-r}^{\tpl} = \llave{\landa_{k-r}(x,y) \, \hpp(x,y) \in \qh_{k+|\tpl|}^{\tpl} : \landa_{k-r} \in  \qh_{k-r}^{\tpl} }.
\eeq
Consider a complement $\Deltap_{k+|\tpl|}$ to  $\hpp\qh_{k-r}^{\tpl}$, i.e.,
$\qh_{k+|\tpl|}^{\tpl}= \Deltap_{k+|\tpl|} \oplus \hpp \qh_{k-r}^{\tpl}$.
Let us introduce the following subspaces of polynomials in two variables (observe that they depend on $\hpp$):
\begin{itemize}
\item $\CCp_k  :=  \llave{ \X_{\g_{k+|\tpl|}} \in \QH_k^{\tpl} : \g_{k+|\tpl|} \in \Deltap_{k+|\tpl|}}.$
\item $\DDp_k := \llave{ \mud_k\,\FZERO \in \QH_k^{\tpl} : \mud_k\in\qh_{k}^{\tpl} }.$
\item $\FFp_k := \llave{ \landa_{k-r}\,\X_{\hpp} \in \QH_k^{\tpl} : \landa_{k-r}\in\qh_{k-r}^{\tpl} }.$
\end{itemize}
A new splitting for planar QH vector fields is presented in the next result.

\begin{pro}\label{descQ}
Let us assume that  $ \hpp \in \qh_{r+|\tpl|}^{\tpl}$ with $\hpp \not\equiv 0$.
Then, $$\QH_k^\tpl = \CCp_k \oplus \DDp_k \oplus \FFp_k.$$
Moreover, each $\P_k \in \QH_k^{\tpl}$  can be uniquely written as
\beq \label{x-div-f}
\P_k=\X_{\g_{k+|\tpl|}} + \mud_k \, \FZERO + \landa_{k-r}\X_\hpp,
\eeq
where
\bea
\nonumber
\g_{k+|\tpl|} &=&\fracp{1}{k+|\tpl|} \Proy_{\Deltap_{k+|\tpl|}}(\FZERO \wedge\P_k) \in \Deltap_{k+|\tpl|},\\
\label{glandamu2}
\mud_k&=&\fracp{1}{k+|\tpl|}\parent{\diverg(\P_k)-\nablapl \landa_{k-r}\cdot\X_\hpp}\in\qh_k^{\tpl}, \\
\nonumber
\landa_{k-r} &=& \fracp{1}{(r+|\tpl|)\hpp} \Proy_{\hpp \qh_{k-r}^{\tpl}}(\FZERO \wedge \P_k)\in\qh_{k-r}^{\tpl}.
\eea
\end{pro}

\begin{proof}
We first show that $\QH_k^\tpl=\CCp_k + \DDp_k + \FFp_k$.
Obviously $\CCp_k + \DDp_k + \FFp_k  \subset \QH_k$.
To prove the converse inclusion, let us to consider $\P_k \in \QH_k^\tpl$.
 From Proposition \ref{con-dis}, we can write
 $\P_k=\X_{\hbar_{k+|\tpl|}} + \mud_k \, \FZERO$ with  $\hbar_{k+|\tpl|}\in \qh_{k+|\tpl|}^\tpl$, $\mud_k \in \qh_k^\tpl$.
 Since $\qh_{k+|\tpl|}^\tpl= \Deltap_{k+|\tpl|} \oplus \hpp \qh_{k-r}^{\tpl}$, we can also write 
 $\hbar_{k+|\tpl|} = \g_{k+|\tpl|} + \landa_{k-r} \hpp$, for some $\g_{k+|\tpl|} \in \Deltap_{k+|\tpl|}$, $\landa_{k-r}\in \qh_{k-r}^\tpl$.
 Then,  we obtain that
 $\P_k =\X_{\g_{k+|\tpl|}} + \X_{\landa_{k-r} \hpp} + \mud_k \, \FZERO$.

Moreover, from Lemma \ref{Oper}\textbf{(b)}, we have
$\landa_{k-r}\X_{\hpp}=\X_{\frac{r+|\tpl|}{k+|\tpl|} \landa_{k-r} \hpp} + \frac{1}{k+|\tpl|}(\nablapl\landa_{k-r} \cdot \X_\hpp)\, \FZERO$,
which implies
$\X_{\landa_{k-r} \hpp} = \frac{k+|\tpl|}{r+|\tpl|} \landa_{k-r}\X_{\hpp} - \frac{1}{r+|\tpl|}(\nablapl \landa_{k-r}\cdot \X_\hpp)\FZERO$.
Consequently:
$$
\P_k  = \X_{\g_{k+|\tpl|}} + \parent{\mud_k - \fracp{1}{r+|\tpl|}(\nablapl \landa_{k-r}\cdot \X_{\hpp})}\FZERO +
\fracp{k+|\tpl|}{r+|\tpl|}\landa_{k-r}\X_\hpp.
$$
To show that the above sum of subspaces is a direct sum, we note that :
\begin{itemize}
\item
$\DDp_k \cap \FFp_k  = \{\zero \}$.
Namely, if $\P_k \in\DDp_k \cap \FFp_k$, then there exist
$\landa_{k-r} \in \qh_{k-r}^\tpl$, $\mud_k \in \qh_{k}^\tpl$ such that $\P_k =\landa_{k-r}\X_\hpp= \mud_k \, \FZERO$.
Hence,
$0=(\mud_k \, \FZERO)\wedge\FZERO= \landa_{k-r}\X_\hpp\wedge\FZERO=(r+|\tpl|)\landa_{k-r}\hpp$.
As $\hpp \not\equiv 0$, we have $\landa_{k-r}=0$ and consequently $\P_k=\zero$.
\item
$\parent{\DDp_k + \FFp_k}\cap\CCp_k = \{\zero \}$.
Namely, if $\P_k \in \parent{\DDp_k + \FFp_k} \cap \CCp_k$, 
then 
$\P_k = \mud_k \, \FZERO+\landa_{k-r}\X_\hpp=\X_{\g_{k+|\tpl|}}$
for some $\landa_{k-r} \in \qh_{k-r}^\tpl$, $\mud_k \in \qh_{k}^\tpl$ and   $\g_{k+|\tpl|}\in\Deltap_{k+|\tpl|}$.
Therefore:
$$
(k+|\tpl|)\g_{k+|\tpl|}=\FZERO\wedge\X_{\g_{k+|\tpl|}}=\FZERO\wedge\P_k=\FZERO\wedge(\mud_k \, \FZERO+\landa_{k-r}\X_\hpp)
=(r+|\tpl|)\landa_{k-r}\hpp.
$$
Hence, $\g_{k+|\tpl|} \in \Deltap_{k+|\tpl|} \cap \hpp \qh_{k-r}^{\tpl}$.
Consequently $\g_{k+|\tpl|}=0$, and then $\P_k=\zero.$
\end{itemize}

To complete the proof, let us  find the expressions for $\g_{k+|\tpl|}$,  $\mud_k$ and $\landa_{k-r}$ given in (\ref{glandamu2}).
Using
$$
\FZERO \wedge \P_k = \FZERO \wedge (\X_{\g_{k+|\tpl|}} + \mud_k\FZERO + \landa_{k-r}\X_\hpp)
= (k+|\tpl|)\g_{k+|\tpl|} + (r+|\tpl|)\landa_{k-r}\hpp,
$$
we easily find the expressions for $\g_{k+|\tpl|}$  and $\landa_{k-r}$.

Moreover, from (\ref{x-div-f})  we have
$\diverg(\P_k)=(k+|\tpl|)\mud_k+ \nablapl\landa_{k-r} \cdot \X_\hpp+\landa_{k-r}\diverg(\X_\hpp)$, and then
$\mud_k=\fracp{1}{k+|\tpl|}\parent{\diverg(\P_k)-\nablapl \landa_{k-r} \cdot \X_\hpp}$.
\end{proof}

\subsection{{
A splitting  for quasi-homogeneous three-dimensional vector fields}}
\label{DecomposTri}

The splitting for planar QH vector fields presented above can be generalized for three-dimensional
QH vector fields. To this end, some definitions are required.

First, we define the  set:
$$
\hpp \qh_{k-r}^{\ttr} = \llave{ \landa_{k-r}(x,y,z)\,\hpp(x,y) \in \qh_{k+|\tpl|}^{\ttr} : \landa_{k-r} \in \qh_{k-r}^{\ttr} },
$$
(notice that it depends on the Hamiltonian function  $\hpp\in \qh_{r+|\tpl|}^{\tpl}$ given in  (\ref{ppalpart}))
 and consider   a complement $\Deltat_{k+|\tpl|}$ to $\hpp \qh_{k-r}^{\ttr}$,
i.e.,  $\qh_{k+|\tpl|}^{\ttr}= \Deltat_{k+|\tpl|} \oplus \hpp \qh_{k-r}^{\ttr}$.

For functions depending on  three variables $h(x,y,z)$,  we  denote
$\X_{h} =\parent{-\frac{\partial h(x,y,z)}{\partial y}\, , \, \frac{\partial h(x,y,z)}{\partial x}}^T.$

Also, we consider the following subspaces  of polynomials in three variables
\begin{itemize}
\item
$\CCt_k = \llave{ \vedor{ \X_{\g_{k+|\tpl|}} }{ 0 } \in \QH_k^{\ttr}  :
\g_{k+|\tpl|} \in \Deltat_{k+|\tpl|}, \g_{k+|\tpl|}(0,0,z)=0 }$.
\item
$\DDt_k = \llave{ \vedor{ \mud_k \FZERO }{ 0 } \in \QH_k^{\ttr} :  \mud_k \in \qh_{k}^{\ttr} }$.
\item
$\FFt_k = \llave{ \vedor{ \landa_{k-r} \X_{\hpp} }{ 0 }\in \QH_k^{\ttr} :
\landa_{k-r} \in \qh_{k-r}^{\ttr} }$.
\item
$\GGt_k = \llave{ \vedor{ \zero }{ \varsig_{k+t_3} } \in \QH_k^{\ttr} :
\varsig_{k+t_3}  \in \qh_{k+t_3}^{\ttr} }$.
\end{itemize}

{
Next result  introduces the  splitting for three-dimensional vector fields.}
\begin{pro}\label{descQtrid}
Let us assume that  $ \hpp\in \qh_{r+|\tpl|}^{\tpl}$ with $\hpp \not\equiv 0$.
Then:
$$\QH_k^{\ttr} = \CCt_k \oplus \DDt_k \oplus\FFt_k \oplus \GGt_k.$$
\end{pro}

\begin{proof}
We first show that $\QH_k^{\ttr} = \CCt_k + \DDt_k +\FFt_k + \GGt_k$.
Obviously $\CCt_k +\DDt_k +\FFt_k + \GGt_k \subset \QH_k^{\ttr}$.

To prove the converse inclusion, let us to consider $\P_k\in \QH_k^{\ttr}$. We can write it as
$$
\P_k = \vedor{ \G_k }{ 0 } + \vedor{ \zero }{ \varsig_{k+t_3} },
$$
where $\varsig_{k+t_3}\in\qh_{k+t_3}^\ttr$  and $\G_k$ depends on the variables $x,y,z$.
We can  write it as:
$$
\G_k(x,y,z)={\sum_{l=0}^{\floor{k/t_3}}}z^l {\Gh}_{k-lt_3}(x,y),
$$
where
$\floor{\, \cdot\, }$ is the \emph{floor} function and ${\Gh}_{k-l t_3}\in \QH_{k-l t_3}^{\tpl}$.
Using Proposition \ref{descQ}, we can write
$$
\G_k(x,y,z)=\sum_{l=0}^{\floor{k/t_3}} z^l \parent{ \X_{\g_{k-lt_3+|\tpl|}}(x,y) + \hmud_{k-lt_3}(x,y)\FZERO +
\hlambda_{k-r-lt_3}(x,y)\X_{\hpp}},
$$
 where
 \bean
\g_{k-lt_3+|\tpl|}&=&\fracp{1}{k-lt_3+|\tpl|}\Proy_{\Deltap_{k-lt_3+|\tpl|}}(\FZERO \wedge {\Gh}_{k-lt_3})
\in \Deltap_{k-l t_3+|\tpl|},\\
\hmud_{k-lt_3}&=&\fracp{1}{k-lt_3+|\tpl|}\parent{\diverg({\Gh}_{k-lt_3})-\nablapl  \hlambda_{k-r-lt_3} \cdot \X_{\hpp}}
\in\qh_{k-l t_3}^{\tpl}, \\
\hlambda_{k-r-lt_3}&=&\fracp{1} {(r+|\tpl|)\hpp}\Proy_{\hpp \qh_{k-r-lt_3}^{\tpl}}  (\FZERO\wedge{\Gh}_{k-lt_3})
\in\qh_{k-l t_3-r}^{\tpl}.
\eean
Therefore $\G_k= \X_{\g_{k+|\tpl|}} + \mud_k \FZERO + \landa_{k-r} \X_{\hpp},$
where
\bean
\g_{k+|\tpl|}(x,y,z) &=& {\sum_{l=0}^{\floor{k/t_3}}}z^l \g_{k-lt_3+|\tpl|}(x,y)\in \Deltat_{k+|\tpl|}\subset \qh_{k+|\tpl|}^{\ttr},\\
\mud_k(x,y,z) &=& {\sum_{l=0}^{\floor{k/t_3}}}z^l \hmud_{k-l t_3}(x,y) \in\qh_{k}^{\ttr}, \\
\landa_{k-r}(x,y,z) &=& {\sum_{l=0}^{\floor{k/t_3}}}z^l\hlambda_{k-r-lt_3}(x,y) \in\qh_{k-r}^{\ttr}.
\eean
We notice that $\g_{k+|\tpl|}(0,0,z)=0$.

To complete the proof,  we observe that $\CCt_k \cap \GGt_k =\{\zero\}$, $\DDt_k\cap \GGt_k =\{\zero\}$ and
$\FFt_k \cap \GGt_k =\{\zero\}$ (this can be shown by reasoning as in the proof of Proposition \ref{descQ}).
Moreover, we easily obtain 
$\DDt_k^{\ttr} \cap \FFt_k^{\ttr}  =\{\zero \}$ and  $(\DDt_k^{\ttr} +\FFt_k^{\ttr})\cap \CCt_k^{\ttr}=\{\zero\}$. 
\end{proof}

\begin{rem}
The above proposition warrants that  each $\P_k\in\QH^{\ttr}_{k}$ can be uniquely written as
$$
\P_k = \vedor{ \X_{\g_{k+|\tpl|}} }{ 0 }  + \vedor{ \mud_k\FZERO }{ 0 }  +
\vedor{ \landa_{k-r}\X_\hpp }{ 0 } + \vedor{ \zero }{ \varsig_{k+t_3} } ,$$
where $\g_{k+|\tpl|}\in \Deltat_{k+|\tpl|}\subset \qh_{k+|\tpl|}^{\ttr}$, $\mud_k\in\qh_{k}^{\ttr}$, $\landa_{k-r} \in\qh_{k-r}^{\ttr}$,
$\varsig_{k+t_3} \in\qh_{k+t_3}^{\ttr}$.
\end{rem}

Next lemma  is useful in the study of the  homological operator. Namely, it characterizes the Lie product $[\P_k,\F_r]$
for any  $\P_k\in\QH^{\ttr}_{k}$, by using the splitting of Proposition \ref{descQtrid}.

\begin{lemma}\label{calcompP}
The following properties hold:
\begin{enumerate}
\item
If $\g_{k+|\tpl|}\in \Deltat_{k+|\tpl|}$ verifies $\g_{k+|\tpl|}(0,0,z)=0$,
then
$$
\bracket{ \vedor{ \X_{\g_{k+|\tpl|}} }{ 0 } , \F_r } = \vedor{ \X_{\tg_{r+k+|\tpl|}} }{ 0 } + \vedor{ \tmu_{r+k} \FZERO }{ 0 } 
+ \vedor{ \tlambda_k \X_{\hpp} }{ 0 } - \vedor{ \zero } {\nablapl \fpp \cdot \X_{\g_{k+|\tpl|}}  },
$$
where
\bea
\nonumber
\tg_{r+k+|\tpl|} &=& \Proy_{\Deltat_{r+k+|\tpl|}}\parent{ \nablapl \g_{k+|\tpl|} \cdot \X_{\hpp}
+ \fracp{k+|\tpl|}{r+k+|\ttr|} \fpp \frac{\partial \g_{k+|\tpl|}}{\partial z}  } \in \Deltat_{r+k+|\tpl|},
\\
\label{glandamu}
\tmu_{r+k} &=& \fracp{1}{r+k+|\tpl|}\parent{\nablapl \fpp \cdot \X_{\frac{\partial \g_{k+|\tpl|}}{\partial z}}
-  \nablapl \tlambda_k\cdot\X_{\hpp}} \in\qh_{r+k}^{\ttr},
\\
\nonumber
\tlambda_k &=& \fracp{r+k+|\tpl|}{(r+|\tpl|)\hpp}\Proy_{\hpp \, \qh_k^{\tpl}}\parent{\nablapl \g_{k+|\tpl|}\cdot \X_{\hpp}
+ \fracp{k+|\tpl|}{r+k+|\ttr|} \fpp \frac{\partial \g_{k+|\tpl|}}{\partial z} } \in\qh_{k}^{\ttr}.
\eea

\item
If $\mud_k\in\qh_{k}^{\ttr}$, then
$\bracket{ \vedor{\mud_k \FZERO }{0} , \F_r } =  \vedor{ (\nablatr \mud_k \cdot \F_r )\FZERO } 0  -
\vedor{ r \mud_k \X_\hpp } 0 - \vedor{ \zero }{ (r+t_3)\mud_k \fpp }$.

\item
If  $\landa_{k-r} \in\qh_{k-r}^{\ttr}$, then
$\bracket{ \vedor{ \landa_{k-r} \X_{\hpp} }{ 0 } , \F_r }=
\vedor{ (\nablatr \landa_{k-r} \cdot \F_r) \X_{\hpp} } 0 - \vedor{ \zero }{  \landa_{k-r} \, \nablapl \fpp \cdot \X_{\hpp} }$.

\item
If $\varsig_{k+t_3} \in\qh_{k+t_3}^{\ttr}$, then
$\bracket{\vedor{\zero} {\varsig_{k+t_3} } , \F_r }
= \vedor{\zero}{\nablatr \varsig_{k+t_3}\cdot \F_r }$.
\end{enumerate}
\end{lemma}

\begin{proof}
In the following computations, we will use the generalized Euler identity for QH
scalar functions: if $\mud_l\in \qh_l^{\tpl}$ then $\nablapl \mud_l \cdot \FZERO = l\,\mud_l$.

\begin{enumerate}
\item A direct computation shows that
\bean
 \bracket{ \vedor{ \X_{\g_{k+|\tpl|}} } {0} , \F_r } 
 &=&  D\vedor{ \X_{\g_{k+|\tpl|}} }{ 0 } \vedor{ \X_{\hpp}}{ \fpp } - D \vedor{ \X_{\hpp}}{ \fpp } \vedor{ \X_{\g_{k+|\tpl|}} }{ 0 } \\
&=& \vedor{  \bracket{\X_{\g_{k+|\tpl|}},\X_{\hpp}} + \fpp \, \X_{\frac{\partial \g_{k+|\tpl|}}{\partial z}}  }{ 0 }
- \vedor{ \zero  }{ \nablapl \fpp \cdot \X_{\g_{k+|\tpl|}} }.\eean

Let us denote $\tP_k:=\bracket{\X_{\g_{k+|\tpl|}},\X_{\hpp}}+ \fpp \, \X_{\frac{\partial \g_{k+|\tpl|}}{\partial z}}$.
From Lemma \ref{Oper}, we can write
$$
\tP_k = \X_{\nablapl \g_{k+|\tpl|} \cdot \X_{\hpp} + \fracp{k+|\tpl|}{r+k+|\ttr|} \fpp \frac{\partial \g_{k+|\tpl|}}{\partial z}}
+ \fracp{1}{r+k+|\ttr|} \parent{ \nablapl \fpp \cdot \X_{\frac{\partial \g_{k+|\tpl|}}{\partial z}} } \FZERO.
$$
To complete the proof, it is enough to  use the splitting  (\ref{descQ}) in order to write
$\tP_k= \X_{\tg_{r+k+|\tpl|}} + \tmu_{r+k}\FZERO + \tlambda_k \X_{\hpp}$.
To obtain the expressions for $\tg_{r+k+|\tpl|}$, $\tmu_{r+k}$ and $\tlambda_k$
given in (\ref{glandamu}), we observe that
$$
\FZERO \wedge \tP_k = (r+k+|\tpl|) \parent{ \nablapl \g_{k+|\tpl|} \cdot \X_{\hpp}
+ \fracp{k+|\tpl|}{r+k+|\ttr|} \fpp \frac{\partial \g_{k+|\tpl|}}{\partial z} }, \ \
\diverg \left(\tP_k\right) = \nablapl \fpp \cdot \X_{\frac{\partial \g_{k+|\tpl|}}{\partial z}}.
$$
\item
It is a straightforward computation.
\item
We have
\bean
\bracket{ \vedor{ \landa_{k-r}\X_{\hpp} }{ 0 } , \F_r } 
&=&
\vedor{ [\landa_{k-r} \X_{\hpp},\X_{\hpp}] + \frac{\partial \landa_{k-r}}{\partial z} \fpp \X_{\hpp} }{ - \nablapl \fpp \cdot \landa_{k-r} \X_{\hpp} }
\\ &=& 
\vedor{ (\nablatr \landa_{k-r} \cdot \F_r)  \X_{\hpp} }{ 0 }
- \vedor{ \zero }{\landa_{k-r}  \nablapl \fpp \cdot \X_{\hpp} }.\eean
\item
As in item 2, it is a straightforward computation.
\end{enumerate}
\end{proof}

\section{Normal forms under  orbital equivalence for  a class of three-dimensional vector fields}
\label{NormalForm}

In this section, we will analyze the orbital homological operator
\bean
\Lequiv_{r+k} &:& \QH_k^{\ttr}
\times \Cor(\Oplt_{k}) \longrightarrow \QH_{r+k}^{\ttr} \\
&& (\P_k,\mut_k) \longrightarrow \Lequiv_{r+k} (\P_k,\mut_k) = [\P_k,\F_r] - \mut_k\,\F_r.
\eean
From Proposition \ref{descQtrid}, this operator can be written as
\bean
\Lequiv_{r+k} &:&
\CCt_k \oplus \DDt_k \oplus \FFt_k \oplus \GGt_k
\times
\Cor(\Oplt_{k}) \longrightarrow \CCt_{r+k} \oplus \DDt_{r+k}
\oplus \FFt_{r+k} \oplus \GGt_{r+k}, \eean
where
\bean
&&\Lequiv_{r+k} \parent{ 
\vedor{ \X_{\g_{k+|\ttr|}} }{ 0 } 
+ \vedor{ \mud_k\FZERO }{ 0  }
+ \vedor{ \landa_{k-r}\X_\hpp }{ 0  }
+ \vedor{ \zero }{ \varsig_{k+t_3} } ,\mut_k 
}
\\
&&\hspace{0.2cm}
=
\bracket{ \vedor{ \X_{\g_{k+|\ttr|}} }{ 0 } , \F_r }
+ \bracket{ \vedor{ \mud_k\FZERO }{ 0 } , \F_r }
+ \bracket{ \vedor{ \landa_{k-r}\X_\hpp }{ 0 } , \F_r }
+ \bracket{\vedor{ \zero }{ \varsig_{k+t_3} } , \F_r } 
- \mut_k\,\F_r.
\eean

As $\mut_k\,\F_r = \vedor{ \mut_k\X_{\hpp} }{ 0 } + \vedor{ \zero }{\mut_k \, \fpp }$,
from Lemma \ref{calcompP}  we obtain that the  matrix associated with  the orbital homological operator $\Lequiv_{r+k}$ is:
$$\arraycolsep=1.2pt 
\begin{array}{|c|c|c|c|c"c|}
\hline
\fantasm \vedor{ \X_{\tg_{r+k+|\tpl|}} }{ 0 }\fantasm 
& 
\vedor{ \zero }{0 } 
&
\vedor{ \zero }{0 } 
&
\vedor{ \zero }{0 } 
& 
\vedor{ \zero }{0 } & \CCt_{r+k}
\\
\hline
\fantasm \vedor{ \tmu_{r+k}\FZERO }{ 0 }\fantasm 
&
\vedor{ (\nablatr \mud_k \cdot \F_r )\FZERO }{ 0 } 
&
\vedor{ \zero }{0 } 
& 
\vedor{ \zero }{0 } 
& 
\vedor{ \zero }{0 } 
& 
\DDt_{r+k}
\\
\hline
\fantasm \vedor{ \tlambda_k \X_{\hpp} }{ 0 }\fantasm 
&
\vedor{ r \mud_k \X_\hpp }{ 0 } 
&
\vedor{ (\nablatr \landa_{k-r} \cdot \F_r) \X_{\hpp} }{ 0 } 
&
\vedor{ \mut_k\X_{\hpp} }{ 0 } 
& 
\vedor{ \zero }{0 } 
& 
\FFt_{r+k}
\\
\hline
\fantasm \vedor{ \zero }{ -\nablapl \fpp \cdot \X_{\g_{k+|\ttr|}} } \fantasm
&
\vedor{ \zero }{ -(r+t_3)\mud_k \fpp } 
&
\vedor{ \zero }{ -\landa_{k-r} \, \nablapl \fpp \cdot \X_{\hpp} } 
&
\vedor{ \zero }{ \mut_k \, \fpp  } 
&
\vedor{ \zero }{ \nablatr \varsig_{k+t_3}\cdot \F_r } 
& 
\GGt_{r+k}
\\
\hlinewd{2pt}
\fantasm \vedor{ \X_{\g_{k+|\ttr|}} }{ 0 } \in \CCt_k\fantasm 
&
\vedor{ \mud_k\FZERO }{ 0 }\in\DDt_k
&
\vedor{ \landa_{k-r}\X_{\hpp} }{ 0 }\in \FFt_k
& 
\mut_k\in\Cor(\Oplt_k) &
\ \
\vedor{ \zero }{  \varsig_{k+t_3} }\in \GGt_k
\ \
&
\\
\hline
\end{array}
$$
where $\tg_{r+k+|\tpl|}$, $\tmu_{r+k}$ and $\tlambda_k$ are given in (\ref{glandamu}).

Using  the Lie derivative operator, the matrix for $\Lequiv_{r+k}$ can also be written as
$$\arraycolsep=1.2pt
\begin{array}{|c|c|c|c|c"c|}
\hline
\fantasm \vedor{ \X_{\tg_{r+k+|\tpl|}} }{ 0 }\fantasm 
& 
\vedor{ \zero }{0 } 
& 
\vedor{ \zero }{0 } 
& 
\vedor{ \zero }{0 } 
& 
\vedor{ \zero }{0 } 
& 
\CCt_{r+k}
\\
\hline
\fantasm \vedor{ \tmu_{r+k}\, \FZERO }{ 0 } \fantasm
&
\vedor{ \Oplt_{r+k}(\mud_k) \,\FZERO }{ 0 } 
& 
\vedor{ \zero }{0 } 
& 
\vedor{ \zero }{0 } 
& 
\vedor{ \zero }{0 } 
& 
\DDt_{r+k}
\\
\hline
\fantasm \vedor{ \tlambda_k\X_{\hpp} }{ 0 }\fantasm 
&
\vedor{ r \mud_k \X_\hpp }{ 0 } 
& 
\vedor{ \Oplt_{r+k}(\landa_{k-r}) \, \X_{\hpp} }{ 0  } 
&
\vedor{ \mut_k\X_{\hpp} }{ 0  } 
& 
\vedor{ \zero }{0 }  
& 
\FFt_{r+k}
\\
\hline
\fantasm \vedor{ \zero }{ -\nablapl \fpp \cdot \X_{\g_{k+|\tpl|}} } \fantasm
&
\vedor{ \zero }{ -(r+t_3)\mud_k  \fpp  } 
&
\vedor{ \zero }{ -\landa_{k-r} \,\nablapl \fpp \cdot \X_{\hpp} } 
&
\vedor{ \zero }{\mut_k \, \fpp  } 
&
\vedor{ \zero }{ \Oplt_{r+k+t_3}(\varsig_{k+t_3}) } 
& 
\GGt_{r+k}
\\
\hlinewd{2pt}
\fantasm \vedor{ \X_{\g_{k+|\tpl|}} }{ 0 } \in \CCt_k \fantasm
& 
\vedor{ \mud_k\FZERO }{ 0 }\in\DDt_k
&
\vedor{ \landa_{k-r}\X_{\hpp} }{ 0 }\in \FFt_k
& 
\mut_k\in\Cor(\Oplt_k) 
& 
\ \ \vedor{ \zero }{ \varsig_{k+t_3} }\in \GGt_k \ \
&  
\\
\hline
\end{array}
$$
where
\bean
\tg_{r+k+|\tpl|}&=&\Opltc_{r+k+|\tpl|,A_0}\parent{\g_{k+|\tpl|}},\\
\tmu_{r+k}&=&\fracp{1}{r+k+|\tpl|}\parent{\nablapl \fpp \cdot \X_{\frac{\partial \g_{k+|\ttr|}}{\partial z}}
-\nablapl \tlambda_k \cdot \X_{\hpp}}
, \\
\tlambda_k&=&\fracp{r+k+|\tpl|}{(r+|\tpl|)\hpp} \Proy_{\hpp \, \qh_k^{\tpl}}\parent{\Oplt_{r+k+|\tpl|,A_0}\parent{\g_{k+|\tpl|}}}
.
\eean
Here, we have denoted  $A_0=\frac{k+|\tpl|}{r+k+|\ttr|}$ and  
 introduced the  modified Lie derivative operator:
\bean
 \Opltc_{k,A} &:& \triangle_{k-r} \longrightarrow \triangle_k \nonumber
\\ &&
\muu_{k-r} \longrightarrow \Proy_{\triangle_k}\parent{\nablapl \muu_{k-r}\cdot\X_{\hpp} + A\, \fpp \, \frac{\partial \muu_{k-r}}{\partial z}}.
\eean

\begin{rem} \label{remark4.10}
If we define the linear operator
\bean
 \Oplt_{k,A} &:& \qh_{k-r}^\t \longrightarrow \qh_k^\t \nonumber
\\ &&
\muu_{k-r} \longrightarrow \nablapl \muu_{k-r}\cdot\X_{\hpp} + A\, \fpp \, \frac{\partial \muu_{k-r}}{\partial z},
\eean
then $\Opltc_{k,A}=\Proy_{\triangle_{k}}\parent{\Oplt_{k,A} |_{\triangle_{k-r}}}$.
Moreover, we notice that  $\Oplt_{k,A}$ agrees with  the Lie derivative operator associated with the vector field
$\F_r + \vedor{ \zero }{ (A-1) \fpp }$.
Taking $A=1$, we obtain the operator  $\Oplt_k$  defined in (\ref{Opele}).
\end{rem}

In the case $\Ker\parent{\Opltc_{r+k+|\tpl|,A_0}}=\llave{0}$,  the upper left block diagonal
of the above matrix has full rank. This allows to  obtain the expression for a complement  $\Cor(\Lequiv_{r+k})$  presented in the next result.
\begin{pro}\label{CorOpHomlogico}
Let us consider $\F_r$ given in (\ref{ppalpart}) and assume that  $\hpp \not\equiv 0$. 
If $\Ker\parent{\Opltc_{r+k+|\tpl|,A_0}}=\llave{0}$ for all $k\in\N$, then
a complement to the range of the orbital homological operator $\Lequiv_{r+k}$ is
$$
\Cor(\Lequiv_{r+k}) = 
\vedor{ \X_{\Cor\parent{\Opltc_{r+k+|\tpl|,A_0}}} }{ 0 } \oplus
\vedor{ \Cor (\Oplt_{r+k})\, \FZERO }{ 0 } \oplus 
\vedor{ \zero }{ \Cor(\Oplt_{r+k+t_3}) },$$
where $A_0=\frac{k+|\tpl|}{r+k+|\ttr|}$  and  $\Cor (\Opltc_{r+k,A_0})$ is a complement to the range of $\Opltc_{r+k,A_0}$ in $ \triangle_{r+k}$.
\end{pro}

Next theorem provides formal normal forms, under conjugation and  orbital equivalence,
for system (\ref{sistprinspecifically}).
We notice that the analysis in the case of  conjugation can be performed by  putting  $\mut_k \equiv 0$ in the above matrix.
The normal form under conjugation obtained in the next theorem has been presented in \cite{AlgCL}.

\begin{theorem}\label{NormaFormFinal-} 
If $\Ker\parent{\Opltc_{r+k+|\tpl|,A_0}}=\llave{0}$ for all $k\in\N$, then
a formal normal form under  conjugation for system (\ref{sistprinspecifically}) is $\dot \x=\G(\x)$, where 
$$
\G = \F_r + \sum_{k\geq 1} \parent{
\vedor{ \X_{\g_{r+k+|\tpl|}} }{ 0 }  +  \vedor{ \mud_{r+k}\, \FZERO }{ 0 }  + \vedor{ \landa_k\, \X_{\hpp} }{ 0 } 
+ \vedor{ \zero }{ \varsig_{r+k+t_3} }
},
$$
with $\g_{r+k+|\tpl|}\in \Cor\parent{\Opltc_{r+k+|\tpl|,A_0}}$, $A_0=\frac{k+|\tpl|}{r+k+|\ttr|}$,
$\mud_{r+k} \in \Cor (\Oplt_{r+k})$,  $\landa_k \in \Cor (\Oplt_{k})$ and
$\varsig_{r+k+t_3} \in \Cor ( \Oplt_{r+k+t_3})$.

Moreover, a formal normal form under orbital equivalence for system (\ref{sistprinspecifically}) is $\dot \x=\G(\x)$, where 
\beq\label{ONF}
\G = \F_r + \sum_{k\geq 1} \parent{
\vedor{ \X_{\g_{r+k+|\tpl|}} }{ 0 } 
+ \vedor{ \mud_{r+k}\, \FZERO }{ 0 } 
+ \vedor{ \zero }{ \varsig_{r+k+t_3}  }
},
\eeq
with $\g_{r+k+|\tpl|}\in \Cor\parent{\Opltc_{r+k+|\tpl|,A_0}}$, $\mud_{r+k} \in \Cor (\Oplt_{r+k})$ and
$\varsig_{r+k+t_3} \in \Cor ( \Oplt_{r+k+t_3})$.
\end{theorem}

\subsection{The operator $\Oplt_{r+k,A}$}
\label{cor}

This subsection is devoted to the study of the linear operator $\Oplt_ {r + k, A}$ for $A\neq0$.
Recall that $\Oplt_ {r + k, 1}=\Oplt_{r+k}$ is the Lie derivative with respect to the principal part of the vector field
$\F_r$ given in (\ref{sistprinspecifically}).

Let us consider $\mud_k \in \qh_k^{\ttr}$ and  write it as
$$\mud_k(x,y,z)= \sum_{l=0}^{\floor{k/t_3}} z^l \, \muh_{k-l t_3}(x,y),$$
with $\muh_{k-l t_3}\in\qh_{k-l t_3}^\tpl$.
Computing the quotient and the remainder $k_1$, $k_2$ 
of the division of integers $k\div t_{3}$, i.e.:
\beq\label{ka}
k=k_1t_3+k_2, \textrm{ with } 0\leq k_2 < t_3,
\eeq
we can write
$\mud_k= \sum_{l=0}^{k_1} z^l \, \muh_{k-lt_3} = \sum_{l=0}^{k_1} z^l \, \muh_{(k_1-l)t_3+k_2}$.

On the other hand, the elements of  the subspace  $\Ker(\Oplp_{r+k})$  are the $(r+k)$-degree QH polynomial first integrals of 
the vector field $\X_\hpp$, which  are of the form $C\,\hpp^\alpha$, with $C,\alpha\in\R$.
Hence, we have $\alpha=\frac k{r+|\tpl|}$  and
\beq\label{nucleo}
\Ker(\Oplp_{r+k})=\spande\llave{\hpp^{\frac k{r+|\tpl|}}}, \textrm{ if } k \textrm{ is a multiple of } r+|\tpl|.
\eeq
For $\alpha\neq\frac k{r+|\tpl|}$, we have $\Ker(\Oplp_{r+k})= \llave{\zero}$.

Consider a complement   $\spacE_k$ to $\Ker(\Oplp_{r+k})$ in $\qh_k^{\tpl}$, i.e.,
$\qh_k^{\tpl} = \spacE_k \oplus \Ker(\Oplp_{r+k})$.
If we  define  $\spacEt_l=z^l\, \spacE_{k-l t_3} $,
$\spacKt_l=z^l\, \Ker(\Oplp_{r+k-l t_3})$, 
then:
\bean
\qh_k^{\ttr} &=&z^{k_1}  \,\qh_{k-{k_1}t_3}^{\tpl}\oplus \cdots \oplus z\, \qh_{k- t_3}^{\tpl} \oplus \qh_{k}^{\tpl}\\
&=&
z^{k_1} \parent{ \spacE_{k-{k_1}t_3} \oplus \Ker(\Oplp_{r+k-{k_1} t_3}) }\oplus \cdots \oplus z \parent{ \spacE_{k- t_3} \oplus \Ker(\Oplp_{r+k- t_3}) } \oplus\parent{ \spacE_{k} \oplus \Ker(\Oplp_{r+k}) }\\
&=&\parent{\spacEt_{k_1} \oplus \spacKt_{k_1} }\oplus \cdots
\oplus\parent{\spacEt_1 \oplus \spacKt_1 }\oplus
 \parent{\spacEt_0 \oplus \spacKt_0 }\\
 &=&\spacEt_{k_1} \oplus\parent{\spacEt_{k_1-1}\oplus \spacKt_{k_1} }\oplus \cdots
\oplus\parent{\spacEt_1 \oplus \spacKt_2 }\oplus
 \parent{\spacEt_0 \oplus \spacKt_1 }\oplus \spacKt_0.
\eean
Given $\muh_j \in \qh_j^{\tpl}$, we can  write it as $\muh_{j}(x,y) = \muh_{j}^{(1)}(x,y) + \muh_{j}^{(2)}(x,y)$, with
$\muh_{j}^{(1)} \in \spacE_j$, $\muh_{j}^{(2)} \in \Ker(\Oplp_{r+j})$.
Therefore:
\bea \nonumber
\mud_k&=& \sum_{l=0}^{k_1} z^{k_1-l} \, \muh_{lt_3+k_2} = \sum_{l=0}^{k_1} z^{k_1-l}\parent{ \muh_{lt_3+k_2}^{(1)}+ \muh_{lt_3+k_2}^{(2)}}
\\ 
\label{mu_k}
&=&z^{k_1}\muh_{k_2}^{(1)}+\sum_{l=1}^{k_1} \parent{ z^{k_1-l} \muh_{lt_3+k_2}^{(1)}+ z^{k_1-(l-1)}\muh_{(l-1)t_3+k_2}^{(2)} } + \muh_{k}^{(2)}.
\eea
After some computations, we obtain
\bea \nonumber
\Oplt_{r+k,A}(\mud_k) &=&
\nablapl \mud_k \cdot \X_{\hpp} + A \, \fpp \, \frac{\partial \mud_k}{\partial z}\\ \nonumber
&=& z^{k_1} \nablapl \muh_{k_2}^{(1)}\cdot\X_{\hpp}
+\sum_{l=1}^{k_1} z^{k_1-l} \parent{ \nablapl \muh_{lt_3+k_2}^{(1)} \cdot\X_{\hpp}+(k_1-l+1)  A \fpp  \muh_{(l-1)t_3+k_2}^{(2)}  } \\ 
\label{ele de mu_k}
&& + A \fpp \sum_{l=1}^{k_1-1} (k_1-l) z^{k_1-l-1} \muh_{lt_3+k_2}^{(1)} .
\eea
Using this equality, we easily obtain that 
$$
\Proy_{z^{k_1}\qh_{r+k_2}^\tpl}\parent{\Oplt_{r+k,A}\parent{z^{k_1}\muh_{k_2}^{(1)}}}  =  z^{k_1}\,\nablapl \muh_{k_2}^{(1)} \cdot\X_{\hpp}=z^{k_1}\,\Oplp_{r+k_2}\parent{\muh_{k_2}^{(1)}}, 
$$
and
\bean
&&
\Proy_{z^{k_1-l}\qh_{r+k_2+lt_3}^\tpl}\parent{\Oplt_{r+k,A}\parent{z^{k_1-l}\muh_{k_2+l t_3}^{(1)}+z^{k_1-(l-1)} \muh_{k_2+(l-1)t_3}^{(2)}}} 
=  \\&&
\hspace{1cm}= z^{k_1-l} \nablapl \muh_{k_2+l t_3}^{(1)} \cdot\X_{\hpp}  
+ (k_1-l+1) \, A\, \fpp \,  z^{k_1-l} \muh_{k_2+(l-1)t_3}^{(2)}.
\eean
Then,  the matrix associated with   $\Oplt_{r+k,A}$ is a block lower triangular matrix:
$$\arraycolsep=1.4pt\def\arraystretch{1.5}
\begin{array}{|c|c|c|c|c|c"c|}
\hline
                     0 & 0 & 0 & 0 & 0 &  0 &  z^{k_1+1}\,\qh_{r+k- (k_1+1) t_3 }^{\ttr} \\
                          \hline
                     D_{k_1} & 0 & 0 & 0 & 0 & 0 &  z^{k_1}\, \qh_{r+k-k_1 t_3 }^{\tpl} \\
                     \hline
                     \bullet & D_{k_1-1} & 0 & 0 & 0 & 0 &  z^{k_1-1}\, \qh_{r+k-(k_1-1) t_3 }^{\tpl} \\
                     \hline
                     \vdots & \ddots & \ddots & \ddots & \vdots & \vdots &  \vdots \\
                     \hline
                    0 & \vdots & \ddots & D_1 &0 & 0 &   z^{1}\, \qh_{r+k-t_3 }^{\tpl}\\
                     \hline
                     0 & 0 & \cdots & \bullet & D_0 &  0 & z^0 \, \qh_{r+k}^{\tpl}\\
\hlinewd{2pt}
                    \ \ \ \ \ \ \spacEt_{k_1}\ \ \ \ \ \  & \spacEt_{k_1-1}\oplus \spacKt_{k_1} &
                     \ \ \ \ \cdots\ \ \ \ \ \ \  &\ \ \spacEt_{1}\oplus \spacKt_{2} \  \ & \ \ \spacEt_{0}\oplus \spacKt_{1}\ \
                                                     & \ \ \ \   \spacKt_{0}\ \ \ \  & \\ \hline
\end{array}
$$
where  $D_{k_1}$ is the matrix associated with the linear operator $\Proy_{z^{k_1}\qh_{r+k_2}^\tpl}\parent{\Oplt_{r+k,A}|_{\spacEt_{k_1}}}$ and $D_{k_1-l}$  is the matrix associated with the linear operator  $\Proy_{z^{k_1-l}\qh_{r+k_2+lt_3}^\tpl}\parent{\Oplt_{r+k,A}|_{\spacEt_{k_1-l}\oplus \spacKt_{k_1-l+1}}}$.
 

\begin{lemma}\label{lerangkerl0}
For each $k\in\N$., we have:  $\Range(\Oplp_{r+k}) \cap \fpp \,\Ker\parent{\Oplp_{r+k-t_3}}=\llave{0}$.
\end{lemma}

\begin{proof}
Recall that we are assuming $ \fpp \in\Cor\parent{\Oplp_{r+t_3}}$ (see Proposition \ref{finCorl}). 
Moreover, if  $k-t_3$ is a multiple of $r+|\tpl|$ we have $\Ker\parent{\Oplp_{r+k-t_3}}=\spande\llave{ \hpp^{\frac{k-t_3}{r+|\tpl|}}}$; otherwise we have $\Ker\parent{\Oplp_{r+k-t_3}}=\llave{ 0}$. 
It is enough to use (\ref{Corciclico}) to obtain the result.
\end{proof}

\begin{lemma}\label{leKerlA}
If $k\in\N$ and $A\in\R\setminus\llave{0}$, then:
$$
\Ker\parent{\Oplt_{r+k,A}}=\Ker\parent{\Oplp_{r+k}} =
\spande\llave{\hpp^n}, \textrm{ if } \ k=n(r+|\tpl|)   \textrm{ for some } n\in\N.$$
Otherwise, we have
$\Ker\parent{\Oplt_{r+k,A}}=\Ker\parent{\Oplp_{r+k}} =\{0\}$.
\end{lemma}

\begin{proof} 
It is trivial to show that $\Ker\parent{\Oplp_{r+k}}=\spacKt_{0}\subset\Ker\parent{\Oplt_{r+k,A}}$.
Let us prove the converse inclusion.
Let us consider $\mud_{k}\in\Ker\parent{\Oplt_{r+k,A}}$ and write it as in (\ref{mu_k}).
From (\ref{ele de mu_k}), we obtain 
\bean 
\Oplt_{r+k,A}(\mud_k)
&=&
z^{k_1}\nablapl \muh_{k_2}^{(1)} \cdot\X_{\hpp} \\&& +\sum_{l=1}^{k_1} z^{k_1-l}\left(\nablapl \muh_{lt_3+k_2}^{(1)} \cdot\X_{\hpp}+(k_1-l+1)A \fpp 
 \left(\muh_{(l-1)t_3+k_2}^{(1)}+\muh_{(l-1)t_3+k_2}^{(2)}  \right)\right)=0.
\eean
Then, we obtain $\nablapl \muh_{k_2}^{(1)} \cdot\X_{\hpp}=0$,  and 
\beq \label{ultima}
\nablapl \muh_{lt_3+k_2}^{(1)} \cdot\X_{\hpp}+(k_1-l+1)A \fpp \left(\muh_{(l-1)t_3+k_2}^{(1)}+\muh_{(l-1)t_3+k_2}^{(2)}\right)=0,
\eeq 
for $l=1,\dots, k_1$.
Consequently:
\begin{itemize}
\item
As $\muh_{k_2}^{(1)}\in\Ker(\Oplp_{r+k_2})\cap\spacE_{k_2}=\llave{0}$, we obtain $\muh_{k_2}^{(1)}=0$.
\item
The equality (\ref{ultima}) with  $l=1$ yields $\nablapl \muh_{t_3+k_2}^{(1)} \cdot\X_{\hpp}=-k_1A \fpp \muh_{k_2}^{(2)}$.
From Lemma \ref{lerangkerl0}, we obtain $\muh_{t_3+k_2}^{(1)}=\muh_{k_2}^{(2)}=0$.
\item
From the equality (\ref{ultima}) with $l=2$ we have 
$$\nablapl \muh_{2t_3+k_2}^{(1)} \cdot\X_{\hpp}+(k_1-1)A \fpp \left(\muh_{t_3+k_2}^{(1)}+\muh_{t_3+k_2}^{(2)}\right)=0,$$ 
that is,  $\nablapl \muh_{2t_3+k_2}^{(1)} \cdot\X_{\hpp}=-(k_1-1)A \fpp \muh_{t_3+k_2}^{(2)}$. 
From Lemma \ref{lerangkerl0}, we get $\muh_{2t_3+k_2}^{(1)}=\muh_{t_3+k_2}^{(2)}=0$.
\item
Reasoning analogously, we obtain $\muh_{lt_3+k_2}^{(1)}=\muh_{(l-1)t_3+k_2}^{(2)}=0$ for $l=0,\dots, k_1-1$. 
\end{itemize}
Therefore $\mud_k=\muh_{k_1t_3+k_2}^{(2)}$, and then we have $\Ker\parent{\Oplt_{r+k,A}}\subset\Ker\parent{\Oplp_{r+k}}$.
\end{proof}

From the structure of the above matrix and Lemma \ref{lerangkerl0}, we state the following result:
\begin{pro}\label{procorlrkA}
Cconsider a complement  $\spacV_{r+k}$ to $\Range(\Oplp_{r+k}) \oplus  \fpp \, \Ker\parent{\Oplp_{k-t_3}}$ in $\qh_{r+k}^{\tpl}$ and denote $k_1=\floor{\fracp{k}{t_3}}$.
Then, a complement  to the range of the operator $\Oplt_{r+k,A}$  is
$$
\Cor\parent{\Oplt_{r+k,A}} = z^{k_1+1} \qh_{r+k-(k_1+1)t_3}^{\ttr}\oplus z^{k_1}\Cor(\Oplp_{r+k_2})\oplus \sum_{l=0}^{k_1-1} z^l \,
 {\spacV}_{r+k-l t_3}.
$$
\end{pro}

\begin{rem}
In the case $\fpp\equiv 0$, we have  ${\spacV}_{r+k}= \Cor(\Oplp_{r+k})$.
\end{rem}

\begin{coro}\label{coCorlrkArk} 
If $k\in\N$ and $A\in\R\setminus\llave{0}$, then, 
$$\Cor\parent{\Oplt_{r+k,A}}=\Cor(\Oplt_{r+k}).$$
\end{coro}
\begin{proof} 
For $A\neq 0$,  the complement $\Cor(\Oplt_{r+k,A})$ to the range
of  $\Oplt_{r+k,A}$ does not  depend on $A$ (see  Proposition \ref{procorlrkA}). 
Hence, $\Cor(\Oplt_{r+k,A})=\Cor(\Oplt_{r+k,1})$, which completes the proof.
\end{proof}

As a consequence of Proposition \ref{procorlrkA}, we obtain:
\begin{theorem}\label{tecorlylc} 
Let us consider  $k\in\N$ and $A_0=\frac{k+|\tpl|}{r+k+|\ttr|}$. 
Let us define $k_1=\floor{\fracp{k}{t_3}}$ and denote by  $\spacV_{r+k}$  a complementary subspace to $\Range(\Oplp_{r+k}) \oplus \fpp \, \Ker\parent{\Oplp_{k-t_3}}$ in $\qh_{r+k}^{\tpl}$.
Then:
\begin{description}
\item[(a)] 
A complement  to the range of  operator $\Oplt_{r+k}$ is
$$
\Cor\parent{\Oplt_{r+k}} = z^{k_1+1} \qh_{r+k-(k_1+1)t_3}^{\ttr}\oplus z^{k_1}\Cor(\Oplp_{r+k_2})\oplus \sum_{l=0}^{k_1-1} z^l \,
 {\spactV}_{r+k-l t_3},
$$
\item[(b)] 
A complement to the range of   operator $\Opltc_{r+k,A_0}$ is
$$
\Cor\parent{\Opltc_{r+k,A_0}} = z^{k_1+1} \Deltap_{r+k-(k_1+1)t_3}\oplus \sum_{l=0}^{k_1} z^l \,
 \Cor\parent{\Proy_{\Deltat_{r+k-l t_3}}\Oplp_{r+k-l t_3}|_{\Deltap_{k-lt_3}}}.
$$
\item[(c)] 
$\Ker\parent{\Opltc_{r+k,A_0}} =\llave{0}$.
\end{description}
\end{theorem}

\begin{proof}
\begin{description}
\item[(a)] 
It is enough to use Proposition \ref{procorlrkA} and Corollary \ref{coCorlrkArk}.

\item[(b)] 
As pointed out in Remark \ref{remark4.10}, we have  
$\Opltc_{r+k,A_0}=\Proy_{\triangle_{r+k}}\parent{\Oplt_{r+k,A_0} |_{\triangle_{k}}}$. 
Moreover, as
\bean
\qh_{k}^{\ttr} &=&
\qh_{k}^{\tpl}
\oplus z \qh_{k-t_3}^{\tpl}
\oplus\cdots
\oplus z^{k_1}\qh_{k-{k_1}t_3}^{\tpl} \\
&=&
 \parent{\Deltap_{k}\oplus \hpp \qh_{k-|\tpl|}^{\tpl}} \oplus
z \parent{\Deltap_{k-t_3}\oplus \hpp \qh_{k-t_3-|\tpl|}^{\tpl}} \oplus \cdots
\oplus
z^{k_1} \parent{\Deltap_{k-{k_1}t_3}\oplus \hpp \qh_{k-{k_1}t_3-|\tpl|}^{\tpl}}, \eean
we obtain 
$\Deltat_{k} = \bigoplus_{l=0}^{k_1}z^l \Deltap_{k-lt_3}$. 
Using Proposition \ref{procorlrkA} and Corollary \ref{coCorlrkArk}, we complete the proof.

\item[(c)] 
Let us suppose on the contrary that $\Ker\parent{\Opltc_{r+k,A_0}}\neq \llave{0}$. 
Then,  there exists $\mud_k\in\Deltat_{k}\setminus\llave{0}$ such that 
$\Proy_{\Deltat_{r+k}}\parent{\Oplt_{r+k,A_0}\parent{\mud_k}}=0$. 
As the matrix associated with the linear operator $\Oplt_{r+k,A_0}$ is triangular, 
we obtain that  $\Proy_{\Deltap_{r+k}}\parent{D_{l}}=0$ for some $l$, with $0\leq l\leq k_1$. 
Moreover, as  $\Deltap_{k-lt_3}\subset\spacE_{k-lt_3}$, then there exists 
$\muh_{k-lt_3}\in\Deltap_{k-lt_3}\setminus\llave{0}$ such that
\bean
\Opltc_{r+k,A_0}\parent{z^l \muh_{k-lt_3}} &=& \Proy_{z^l\Deltap_{r+k-lt_3}}\parent{\Oplt_{r+k,A_0}(z^l\muh_{k-lt_3})}
\\ &=& z^l \Proy_{\Deltap_{r+k-lt_3}}\parent{\nablapl \muh_{k-lt_3}\cdot\X_{\hpp}}=0.
\eean
Hence, $\nablapl \muh_{k-lt_3}\cdot\X_{\hpp}=\eta \hpp$, for some  $\eta\in\qh_{k-lt_3-r-|\tpl|}^{\tpl}$. 
Therefore, $\nablapl\hpp \cdot \X_{\muh_{k-lt_3}}=\eta\hpp$, that is, $\hpp$ is an invariant curve for $\X_{\muh_{k-lt_3}}$.
Consequently, $\muh_{k-lt_3}\in\Deltap_{k-lt_3}\cap \spande\llave{ \hpp }=\llave{0}$, which is a  contradiction.
\end{description}
\end{proof}

A normal form for system (\ref{sistprinspecifically}), whose principal part is given in (\ref{ppalpart}), can be easily derived by 
combining  Theorems \ref{NormaFormFinal-}  and \ref{tecorlylc}.

\section{The normal form for the Hopf-zero singularity}
\label{Hopzero}

In this section we apply the above ideas  to obtain a normal form for a Hopf-zero singularity.
As we will see, we  reach simpler normal forms than the classical ones obtained by using Taylor expansions of the vector field (see, e..g., \cite{Guckenheimer83b,Kuz95}). 
Following similar ideas, the normal form under conjugation for a triple zero degeneracy has been obtained in \cite{AlgCL}.

The study of further simplifications in the classical normal form under conjugation for the Hopf-Zero singularity was initiated in \cite{Ushiki84},
where simplified normal forms were obtained up to low order.
Later, the study of such further simplifications up to an arbitrary order was considered in \cite{AlgabaHZ},
where the use of orbital equivalence was also included.

More recently,  some authors have addressed the problem of obtaining the simplest
normal form under orbital equivalence for the Hopf-Zero singularity under generic conditions.
Specifically, in \cite{CWY2003,CWY2005},  the simplest (unique) normal forms under
conjugation and orbital equivalence, for system  (\ref{sisprinc}) is obtained in the nondegenerate case.
On the other hand, the study  of the simplest orbital normal form for free divergence systems with Hopf-zero singularity is carried out in \cite{GM2013}, whereas in \cite{GMS2013} the case corresponding to the presence of an invariant plane is considered.

Let us consider a smooth three-dimensional system having an equilibrium point at the origin, with eigenvalues  
$0$, $\pm i$. 
Using adequate coordinates, the system can be written as 
$$
\vetre{\dot x}{\dot y}{\dot z} = \vetre{-y}{x}{0} +\vetre{f(x,y,z)}{g(x,y,z)}{h(x,y,z)},
$$
where $f,g,h$ denote the nonlinear terms.
We assume that the generic hypothesis $\frac{\partial^2 h }{\partial x^2} \neq 0$ or  $\frac{\partial^2 h }{\partial y^2} \neq 0$ holds.
Selecting the  type $\ttr=(1,1,2)$, it is easy to show that,  after a  suitable 0th-degree transformation and a rescaling,  the above system can be expanded in QH  of type $\ttr$ so that its principal part becomes:
$$
\vetrer{ -y }{ x }{\frac{x^2+y^2}2 } =
\vedor{ \X_{\hpp} }{ \fpp }\in \QH_0^{\ttr},
$$ 
where 
$$\hpp= \fpp=\frac{x^2+y^2}2\in \qh_{2}^{\tpl}.$$
We observe that, in this case,   $\tpl=(1,1)$, $r=0$ and
$\FZERO=\nabla\hpp= \vedot{x}{y} \in \QH_{0}^{\tpl}$.
We also notice that $\fpp\in\Cor(\Oplp_{2})$, a complement to the range of the Lie derivative operator $\Oplp_{2}$ (see Proposition \ref{finCorl}).
In fact, we have $\Cor(\Oplp_{k}) = \spande\llave{ \hpp^{{k}/{2}}}$,  if $k$ is even, or the trivial space $\llave{0}$,   if $k$   is odd (see  \cite{Algaba12}).

Summarizing, the nondegenerate Hopf-zero singularity expanded in QH terms of type $\ttr=(1,1,2)$ can be written as:
\beq\label{hopf-zero-system}
\dot{\x}=\F(\x)=\F_0(\x)+\F_1(\x)+\cdots,
\eeq
where the principal part is
\beq\label{princpartHZ}
\F_0(\x)=\vetrer{ -y }{ x }{ \frac{x^2+y^2}2 } = \vedor{ \X_{\hpp} }{ \hpp }\in \QH_0^{\ttr}.
\eeq

From Proposition \ref{procorlrkA}, as  $r+|\tpl|=2$, we must consider two cases corresponding to 
$k-t_3=k-2$ even or odd.
After some simple computations, it can be shown that:
\begin{itemize}
\item If   $k=2k_1$ is even, then 
$\frac{k-t_3}{r+|\tpl|}=k_1-1$, 
$\Ker\parent{\Oplp_{2(k_1-j)-2}}=\spande\llave{\hpp^{k_1-1-j}}$ 
and a complementary subspace  to $\Range(\Oplt_{2(k_1-j)}) \oplus \hpp \, \Ker\parent{\Oplp_{2(k_1-j)-2}}$ is 
$\spacV_{2(k_1-j)}=\llave{0}$.
Hence:
\beq\label{corlpar}
\Cor(\Oplt_{2k_1}) 
= z^{k_1+1}\qh_{-2}^{\ttr} \oplus z^{k_1}\Cor(\Oplp_{0}) \oplus {\sum_{j = 0}^{k_1-1}} z^{j} \spacV_{2(k_1-j)}
=\spande\llave{z^{k_1}}.
\eeq

\item If  $k=2k_1+1$ is odd, then  
$\frac{k-t_3}{r+|\tpl|}=\frac{2k_1-1}{2}$, 
$\Ker\parent{\Oplp_{2(k_1-j)-1}}=\llave{0}$ 
and a complementary subspace  to $\Range(\Oplt_{2(k_1-j)+1}) \oplus \hpp \, \Ker\parent{\Oplp_{2(k_1-j)-1}}$ is  
$\spacV_{2(k_1-j)+1})=\llave{0}$. Hence:
\beq\label{corlimpar}
\Cor(\Oplt_{2k_1+1})
= z^{k_1+1}\qh_{-1}^{\ttr} \oplus z^{k_1}\Cor(\Oplp_{1}) \oplus  {\sum_{j = 0}^{k_1-1}} z^{j}\spacV_{2(k_1-j)+1}) =\llave{0}.
\eeq
\end{itemize}

Next result can be derived from Theorems  \ref{NormaFormFinal-} and \ref{tecorlylc}:
\begin{theorem}\label{NFHopf-Zero-}
A formal normal form under conjugation for system (\ref{hopf-zero-system}) is given by
\beq\label{FNHZConj}
\dot{\x} = \F_0(\x)  +  \vedor{ G_1(z)\nabla\hpp + G_3(z)\X_{\hpp} }{ G_2(z) }  ,
\eeq
where $G_1(z)= {\sum_{l \geq 1}^{\infty}} a_{l}z^{l}$, $G_2(z)= {\sum_{l\geq 1}^{\infty}} b_{l}z^{l+1}$ and
$G_3(z)= {\sum_{l\geq 1}^{\infty}} c_{l}z^{l}$.

A formal normal form  under orbital equivalence for system (\ref{hopf-zero-system}) is given by
\beq\label{FNHZEquiv}
\dot{\x} =\F_0(\x) +  \vedor{ G_1(z)\nabla\hpp  }{ G_2(z) } ,
\eeq
where $G_1(z)= {\sum_{l \geq 1}^{\infty}} \tilde{a}_{l}z^{l}$ and  $G_2(z)= {\sum_{l \geq 1}^{\infty}} \tilde{b}_{l}z^{l+1}$.
Moreover, $\tilde{a}_{1}=a_1$ and $\tilde{b}_{1}=b_1$.
\end{theorem}

The normal form (\ref{FNHZEquiv}) is used in \cite{AlgHZInteg} to analyze  the integrability problem for the Hopf-zero  singularity.
Also, the relation of the integrability problem with the existence of inverse Jacobi multipliers is stablished
{
in \cite{AlgHZInteg} by means of the  normal form (\ref{FNHZEquiv}).}

\subsection{The parametric normal form  for the Hopf-zero singularity}

In practice, it  is usual to deal with dynamical systems depending on parameters. 
In the analysis of bifurcations of families of vector fields,  the normal form procedure is of great interest.
For instance, one can obtain the normal form for the critical parameter value where the bifurcation takes place and then
analyze the bifurcation behavior of an adequate unfolding of the normal form.

More accurate results can be obtained by using parametric  normal  forms,
because explicit relations between the original parameters and
the unfolding parameters can be derived.
For instance,  one can suspend the system  with the parameter, treating it 
as new state variable and extending to a higher-dimensional system with no parameters, and calculate normal forms of the
extended system. 

Another possibility, which is simpler, consists in  reducing the principal part to an adequate miniversal unfolding 
and later perform a normal form reduction taking into account the parameters.
Despite the importance of parametric normal forms in applications, there have been very few results in the literature; see
\cite{GazMoaz15,YuLeung03,GazSadri16,GaoZhang10}.

The parametric normal form  for the Hopf-zero singularity is presented in the following results.

\begin{pro} \label{ProposMiniversal}
A miniversal unfolding  for the system associated to the principal part (\ref{princpartHZ}):
\bea
\nonumber
\dot x &=& -y, \\
\label{sist1}
\dot y &=& x, \\
\nonumber
\dot z &=& \fracp{1}{2}x^2+\fracp{1}{2}y^2, 
\eea
is given by:
\bea
\nonumber
\dot x &=& -y+\epsilon x, \\
\label{sist2}
\dot y &=& x+\epsilon y, \\
\nonumber
\dot z &=&  \delta z +\fracp{1}{2}x^2+\fracp{1}{2}y^2.
\eea
\end{pro}

\begin{proof}
A parametric versal unfolding for system (\ref{sist1}) is given by:
\bea
\nonumber
 \dot{x} &=&
 -(1+\alpha_2)y+\alpha_1 x, \\
 \label{sist3}
\dot{y} &=&
  (1+\beta_2)x+\beta_1 y,\\
\nonumber
 \dot{z}
 &=&
  \fracp{1}{2}x^2+\fracp{1}{2}y^2+\gamma_1 x+\gamma_2 y+\gamma_3 z+\gamma_4 xy+\gamma_5x^2+\gamma_6y^2.
  \eea

 We will perform a sequence of rescaling and 0th-degree QH transformations.
First, we introduce the variables 
$$ \txx = x / \sqrt{1+\alpha_2} , \  \tyy = y / \sqrt{1+\beta_2} , \ \tzz = z + \dd_1 x + \dd_2 y,  $$
where 
\bean
\dd_1 &=& 
\parent{{(\beta_1-\gamma_3)\gamma_{1}-(1+\beta_2)\gamma_{2}}} / 
\parent{\alpha_1\beta_1-\alpha_1\gamma_3+\alpha_2\beta_2-\beta_1\gamma_3+\gamma_3^2+\alpha_2+\beta_2+1}, \\ 
\dd_2 &=& 
\parent{{(1+\alpha_2)\gamma_{1}+(\alpha_1-\gamma_3)\gamma_{2}}} /
\parent{\alpha_1\beta_1-\alpha_1\gamma_3+\alpha_2\beta_2-\beta_1\gamma_3+\gamma_3^2+\alpha_2+\beta_2+1}.
\eean
Then,  system (\ref{sist3}) is transformed into
\bea
\nonumber
\dot{\txx}  &=& - \ta_1 \tyy + \alpha_1 \txx,\\
\label{sist4} 
\dot{\tyy} &=& \ta_1 \txx + \beta_1 \tyy,\\ 
\nonumber
\dot{\tzz} &=& \parent{\fracp{1}{2}+\tilde\gamma_5}\txx^2 + \parent{\fracp{1}{2}+\tilde\gamma_6} \tyy^2 + \gamma_3 \tzz + \tgam_4 \txx \tyy,
\eea
where
\bean
\ta_1 &=& \sqrt{1+\alpha_2} \sqrt{1+\beta_2} = 1 + {\cal O}(|\alpha_2,\beta_2|),\\
\tgam_4 &=& \gamma_4 \sqrt{1+\alpha_2} \sqrt{1+\beta_2} = \gamma_4 + {\cal O}(|\alpha_2,\beta_2|),\\
\tgam_5 &=& \parent{2\alpha_2\gamma_5+\alpha_2+2\gamma_5} / {2},\\
\tgam_6 &=& \parent{2\beta_2\gamma_6+\beta_2+2\gamma_6} / {2}.
\eean

Next, we perform the linear transformation 
$$ \hxx = \txx + a \tyy, \ \hyy = a \txx + \tyy, \ \hzz = \tzz, $$
where
$$a= - \parent{2\ta_1+\sqrt{4\ta_1^2-(\alpha_1-\beta_1)^2}}/\parent{\alpha_1-\beta_1},$$ 
 (notice that  $1-a^2\neq 0$).
Then,  system (\ref{sist4}) becomes
\bea 
\nonumber
 \dot{\hxx} &=& - \ha_1 \hyy + \heps \hxx,\\ 
\label{sist5}
\dot{\hyy} &=& \ha_1 \hxx + \heps \hyy,\\ 
\nonumber
\dot{\hzz}  &=& \hc_1 \hxx^2 + \hc_2 \hyy^2 + \gamma_3 \hzz + \hgam_4 \hxx \hyy,
 \eea
where
\bean
\ha_1 &=& - \parent{\ta_1 a^2+(\alpha_1-\beta_1)a+\ta_1} / \parent{a^2-1} = \ta_1 + {\cal O}(a),\\
\heps &=& \parent{\beta_1 a^2 - 2 \ta_1 a-\alpha_1} / \parent{a^2-1} = {\cal O}(\alpha_1,\beta_1),\\
\hc_1 &=& \parent{(1+2\tgam_6) a^2-2\tgam_4 a+1+2\tgam_5} / \parent{2(a^2-1)^2} = \parent{1/2 + \tgam_5} + {\cal O}(a),\\
\hc_2 &=& \parent{(1+2\tgam_5) a^2-2\tgam_4 a+1+2\tgam_6} / \parent{2(a^2-1)^2} = \parent{1/2 + \tgam_6}+ {\cal O}(a),\\
\hgam_4 &=& \parent{\tgam_4 a^2-2(\tgam_5+\tgam_6) a-2a+\tgam_4} / {(a^2-1)^2}=\tgam_4+{\cal O}(a).
\eean
Finally, we use a time-reparametrization $dt={dT}/{\ha_1}$, and a near-identity transformation and a rescaling given by
$$ \bxx = \hxx, \ 
\byy = \hyy , \
\bzz = (\hzz + \ee_1 \hxx \hyy + \ee_2 \hxx^2)/c,$$
where
\bean
\ee_1 &=& - ((2 \heps-\gamma_3)\bgam_4 + 2(\hc_1-\hc_2) {\ha_1} ) / ((2\heps-\gamma_3)^2+4 {\ha_1}^2), \\
\ee_2 &=& - ( (2\heps-\gamma_3)(\hc_1-\hc_2) - 2 \hgam_4 {\ha_1} ) / ((2\heps-\hgam_3)^2+4 {\ha_1}^2), \\ 
 c   &=& 2( (2\heps-\gamma_3)^2 \hc_2 +(2\heps-\hgam_3)\hgam_4 {\ha_1} +2(\hc_1+hc_2) {\ha_1}^2) / (\ha_1((2\heps-\gamma_3)^2+4{\ha_1}^2)).
 \eean
Then, system (\ref{sist5}) is transformed into system (\ref{sist2}) with $\epsilon=\heps/{\ha_1}$, $\delta = \gamma_3/{\ha_1}$.
\end{proof}

\begin{theorem}\label{teoHZFN} 
Let us consider the perturbation of the unfolding (\ref{sist2})
given by
\beq\label{SistEst} 
\vetre{\dot{x}}{\dot{y}}{\dot{z}} =
\vetrer{-y+\epsilon x}{x +\epsilon y}{\delta z +\fracp{1}{2}x^2+\fracp{1}{2}y^2}
+\F_1+\F_2+\cdots,
\eeq
where
\bean
\F_1&=&
\vetrer {A_{200}x^2+A_{110}xy+A_{020}y^2+(A'_{001} \epsilon + A''_{001} \delta)z}
{B_{200}x^2+B_{110}xy+B_{020}y^2+(B'_{001} \epsilon + B''_{001} \delta)z}
{C_{300}x^3+C_{210}x^2y+C_{120}xy^2+C_{102}xz+C_{030}y^3+C_{012}yz} + {\cal O}\left(\left|\epsilon,\delta\right|^2\right) \in \QH_1^{\ttr},
\\ 
\F_2&=& 
\vetrer 
{A_{300}x^3+A_{210}x^2y+A_{120}xy^2+A_{102}xz+A_{030}y^3+A_{012}yz}
{B_{300}x^3+B_{210}x^2y+B_{120}xy^2+B_{102}xz+B_{030}y^3+B_{012}yz}
{\begin{array}{c}
C_{400}x^4+C_{310}x^3y+C_{220}x^2y^2+C_{201}x^2z\\+C_{130}xy^3+C_{111}xyz+C_{040}y^4+C_{021}y^2z+C_{002}z^2
\end{array}} + {\cal O}\left(\left|\epsilon,\delta\right|^2\right) \in \QH_2^{\ttr}.
\eean
Then:
\begin{description}
\item[(a)] 
A parametric normal form for this family is given by:
\beq \label{NF_HoZe}
\vetrer
{\dot{x}}
{\dot{y}}
{\dot{z}} =
\vetrer
{-y+\epsilon x}
{x +\epsilon y}
{ \delta z +\fracp{1}{2}x^2+\fracp{1}{2}y^2}
+a_1 z \vetrer xy0 + b_1\vetrer 00{z^2}+c_1 z \vetrer{-y}x0,
\eeq
where:
\bean 
a_1&=&
\fracp 1{2}  A_{101} + \fracp 1{2}  B_{011}   
\\&&
 +  \epsilon
\left(
   \fracp 1{2}  A'_{001}  A_{110} 
   + A'_{001}  B_{020} 
   - \fracp 1{2}  A'_{001}  C_{011} 
- A_{200} B'_{001}  
+ \fracp 1{2}  B'_{001}  C_{101}
 - \fracp 1{2}  B'_{001}  B_{110}
\right) 
 \\&&
 + 
\delta
\left(
- \fracp 14  A_{120} - \fracp 34  A_{300} - \fracp 34  B_{030} - \fracp 14  B_{210} 
 - \fracp 14  A_{020} A_{110} - \fracp 1{2}  A_{020} B_{020} - \fracp 14  A_{110} A_{200} 
\right. 
\\&& 
\left. 
+ \fracp 1{2}  A_{200} B_{200}  + \fracp 14  B_{020} B_{110} + \fracp 14  B_{110} B_{200}
+ \fracp 1{2}  A''_{001} A_{110} 
+  A''_{001}  B_{020} 
- \fracp 1{2}  A''_{001}  C_{011} 
\right. 
\\&&
\left. 
- A_{200}  B''_{001} 
+ \fracp 1{2}   B''_{001}  C_{101} 
- \fracp 1{2} B''_{001}  B_{110}
\right)
+ {\cal O}\left(\left|\epsilon,\delta\right|^2\right),
\eean
\bean
b_1&=&
C_{002} 
 +  \epsilon 
\left( 
 A'_{001}C_{011} 
 - B'_{001} C_{101} 
\right)
 \\&& 
+ \delta
\left(
  - \fracp 14  A_{120} 
  - \fracp 34  A_{300} 
  - \fracp 34  B_{030} 
  - \fracp 14  B_{210} 
  - \fracp 1{2}  C_{021} 
  - \fracp 1{2}  C_{201}
- \fracp 14  A_{020} A_{110} 
   \right. \\ && \left. 
- \fracp 14  A_{110} A_{200} 
+ \fracp 14  B_{020} B_{110} 
+ \fracp 14  B_{110} B_{200} 
- \fracp 1{2}  A_{020} B_{020} 
+ \fracp 1{2}  A_{200} B_{200} 
- \fracp 1{2}  A_{020} C_{011} 
\right. \\ && \left. 
- \fracp 1{2}  A_{200} C_{011}   
+ \fracp 1{2}  B_{020} C_{101} 
+ \fracp 1{2}  B_{200} C_{101}
 +  A''_{001} C_{011} 
 -  B''_{001} C_{101} 
 \right)
+ {\cal O}\left(\left|\epsilon,\delta\right|^2\right),
\eean
\bean
c_1&=&
 - \fracp 1{2}  A_{011} 
 + \fracp 1{2}  B_{101} 
  \\&& 
 + \epsilon
\left( 
- A'_{001}  A_{020} 
 + \fracp 1{2}  A'_{001} B_{110} 
 - \fracp 1{2}  A'_{001} C_{101} 
 + \fracp 1{2}  A_{110} B'_{001}  
 - \fracp 1{2}  B'_{001} C_{011} 
 - B'_{001} B_{200}
\right)
 \\&&
 + \delta
\left( 
\fracp 34  A_{030} 
+ \fracp 14  A_{210} 
- \fracp 14  B_{120} 
- \fracp 34  B_{300} 
+ \fracp 56  A_{020}^2 
+ \fracp 1{12}  A_{110}^2 
+ \fracp 13  A_{200}^2 
+ \fracp 13  B_{020}^2 
  \right. \\ && \left. 
+ \fracp 1{12}  B_{110}^2 
+ \fracp 56  B_{200}^2 
- \fracp 1{12}  A_{020} B_{110} 
- \fracp 1{12}  A_{110} B_{200} 
+ \fracp 56  A_{020} A_{200} 
- \fracp 5{12}  A_{110} B_{020}   
\right. \\ && \left. 
- \fracp 5{12}  A_{200} B_{110} 
+ \fracp 56  B_{020} B_{200}
- A''_{001}  A_{020} 
 + \fracp 1{2} A''_{001}  B_{110} 
 - \fracp 1{2}  A''_{001}  C_{101} 
 + \fracp 1{2}  A_{110} B''_{001}  
   \right. \\ && \left. 
 - \fracp 1{2} B''_{001} C_{011} 
-  B''_{001} B_{200}
\right)
+ {\cal O}\left(\left|\epsilon,\delta\right|^2\right).
\eean

\item[(b)] 
A parametric  orbital normal form for this family is given by:
\beq \label{ONF_HoZe}
\vetrer
{\dot{x}}
{\dot{y}}
{\dot{z}} =
\vetrer
{-y+\epsilon x}
{x +\epsilon y}
{ \delta z +\fracp{1}{2}x^2+\fracp{1}{2}y^2}
+\hat{a}_1 z \vetrer xy0 + \tilde{b}_1\vetrer 00{z^2},
\eeq
where  $\tilde{a}_{1}=a_1$, $\tilde{b}_{1}=b_1$ are given before.
\end{description}
\end{theorem}

\subsection{Computing normal form coefficients}  \label{ComputCoef}

 \def\U{{\mathbf U}}

Here, we include the proof of Theorem \ref{teoHZFN}. We have included it in a separate subsection because 
it is interesting to show how the expressions for the normal form coefficients,  given in the  theorem,  are obtained.

Let us consider the vector field  
$$\F_0=\vetrer{-y+\epsilon x}{x+\epsilon y}{\fracp{x^2}{2}+\fracp{y^2}{2}+\delta z}\in\QH_0^\t,$$
associated with system (\ref{sist2}).
Then,  the vector field corresponding to system (\ref{SistEst}) can be written as 
$\F(\x)=\F_0(\x)+\F_{1}(\x)+\cdots$, where    $\F_k\in\QH_k^\t$. 

According to Section \ref{ONF},  we first reparametrize the time by $\fracp{dt}{d \tau}=1-\mu(\x)$,
with $\mu=\sum_{j\geq 1}\mu_j$, being $\mu_j\in\Cor\parent{\Opl_j}$.
From  (\ref{corlpar}), (\ref{corlimpar}), we have  $\mu_{2m+1}=0$ and $\mu_{2m}=\gamma_{2m}z^m$, for $m\geq 1$.

Next, we perform a near-identity transformation corresponding to a
generator $\U=\sum_{j\geq 1}\U_j$, where $\U_j\in\QH_j^\t$. Then,
the transformed vector field is:
\bean
\G&=&\parent{1-\mu}\F+\bracket{\parent{1-\mu}\F,\U}+\fracp{1}{2!}\bracket{\bracket{\parent{1-\mu}\F,\U},\U}
\\ &&
+\fracp{1}{3!}\bracket{\bracket{\bracket{\parent{1-\mu}\F,\U},\U},\U}
+\fracp{1}{4!}\bracket{\bracket{\bracket{\bracket{\parent{1-\mu}\F,\U},\U},\U},\U}+\cdots.
\eean

We deal with this vector field degree by degree
\begin{itemize}
\item
It is easy to show that the 0th degree QH terms do not change: $\G_0=\F_0$.
\item
The QH term of degree 1 of the transformed vector field is given by
$$
\G_1 = \F_1+\bracket{\F_0,\U_1}=\F_1-\Lequiv_1\parent{\U_1,0}.
$$
To annihilate this term, we choose $\U_1$ as the  unique solution of the homological equation $\Lequiv_1\parent{\U_1,0}=\F_1$.
To compute this solution, we write $\F_1$ as in the statement of Theorem  \ref{teoHZFN}.
Using a computer algebra system, we have obtained  $\U_1=(p_1,q_1,r_1)^T$, where
\bean
p_1 &=&
- \left(   \fracp13 A_{110} + \fracp23 B_{020} + \fracp13 B_{200} \right) x^2
+\left( \fracp23 A_{200}+ \fracp13 B_{110} - \fracp23 A_{020} \right) x y
 \\ &&
+\left( \fracp13 A_{110} - \fracp13 B_{020}  - \fracp23  B_{200} \right) y^2
+{\cal O}\left(\left|\epsilon,\delta\right|\right),
\\
q_1 &=& 
 \left( \fracp23 A_{020} + \fracp13 A_{200} - \fracp13 B_{110}  \right) x^2
+\left( \fracp23 B_{200} - \fracp13 A_{110} - \fracp23 B_{020} \right) x y\\ &&
+\left( \fracp23 A_{200} + \fracp13 B_{110}  + \fracp13 A_{020} \right) y^2
+{\cal O}\left(\left|\epsilon,\delta\right|\right),
\\
r_1 &=& 
-C_{011} x z
+ C_{101} y z
+ \left( \fracp12 C_{101} -  \fracp29 B_{110} - \fracp23 C_{030} - \fracp29 A_{020} - \fracp79 A_{200} - \fracp13 C_{210}  \right) x^3
\\ &&
+ \left( C_{300} + \fracp12 C_{011} - \fracp23 B_{020} -  \fracp13 B_{200} - \fracp13 A_{110} \right)  x^2 y 
\\ &&
+ \left( \fracp12 C_{101}  - \fracp23 A_{200} - \fracp13 B_{110} - \fracp13 A_{020} - C_{030} \right) x y^2
\\ &&
+ \left( \fracp12 C_{011} - \fracp79 B_{020} + \fracp13 C_{120} - \fracp29 B_{200} + \fracp23 C_{300} - \fracp29 A_{110} \right) y^3
+{\cal O}\left(\left|\epsilon,\delta\right|\right).
\eean

\item
The  QH term of degree 2 of the transformed vector field is given by
$$ \G_2 = -\mu_2\F_0+\bracket{\F_0,\U_2}+\P_2=-\Lequiv_2(\U_2,\mu_2)+\P_2,$$
where
$$
\P_2 = \F_2+\bracket{\F_1,\U_1}+\fracp{1}{2!}\bracket{\bracket{\F_0,\U_1},\U_1}.$$
To obtain the expression of $a_1,b_1$  given in Theorem \ref{NFHopf-Zero-}, it is enough to obtain
$$\P_2^c = a_1 z \vetret{x}{y}{0} + b_1 z^2 \vetret{0}{0}{1},$$ 
the projection of  $\P_2$ onto $\Cor\parent{\Lequiv_2}$.
Next, we choose  $\parent{\U_2,\mu_2}$ as the  unique  solution of the linear system $\Lequiv_2\parent{\U_2,\mu_2}=\P_2^r=\P_2-\P_2^c$
(in this way, we get $\G_2=\P_2^c$). 

By using a computer algebra system, we have computed the analytical expressions of $\U_2=(p_2,q_2,r_2)^T $ 
and $\mu_{2} = a_2 z$, but they are too long to be presented here.

\item
The  QH term of degree 3 of the transformed vector field is given by
$$ \G_3 = \bracket{\F_0,\U_3}+\P_3=-\Lequiv_3\parent{\U_3,0}+\P_3,$$
where
\bean
\P_3&=&\F_3+\bracket{\F_1,\U_2}+\bracket{\F_2-\mu_2\F_0,\U_1}+\fracp{1}{2!}\bracket{\bracket{\F_0,\U_1},\U_2}+\fracp{1}{2!}\bracket{\bracket{\F_0,\U_2},\U_1}\\
&&+\fracp{1}{2!}\bracket{\bracket{\F_1,\U_1},\U_1}+\fracp{1}{3!}\bracket{\bracket{\bracket{\F_0,\U_1},\U_1},\U_1}.
\eean
Again, it is possible to annihilate this term by choosing $\U_3$ as the  unique solution of the homological equation   $\Lequiv_3\parent{\U_3,0}=\P_3$ (see  Proposition \ref{CorOpHomlogico}).
\end{itemize}

\subsection{A case study: the Hopf-zero bifurcation in the Fitzhugh-Nagumo system}

In this last subsection, we apply the results to the three-dimensional  
FitzHugh-Nagumo differential system:
\bea
\nonumber
\dot x &=& z, \\
\label{FitzNag}
\dot y &=& b ( x - d y), \\
\nonumber
\dot z &=& f(x)+y+cz,
\eea
where $f(x)=x(x-1)(x-a)$.
This system appears when looking for traveling waves in the  FitzHugh-Nagumo partial differential system:
$u_t=u_{xx} - f(u)-v, v_t=\delta(u-\gamma v)$, that is a simple model describing the excitation of neural membranes and the propagation of nerve impulses along an axon (see \cite{Euzebio15} and references therein).
Among other situations, system (\ref{FitzNag}) exhibits a  Hopf-zero singularity for the equilibrium at the origin for the  parameters values
$$a^*=-1/d, \ c^*=bd, \ d(1-b^2d^3)>0.$$
For these   critical  values, the linearization matrix 
$$A=\parent{ \begin{array}{ccc} 0 & 0 & 1 \\  b & - b d & 0 \\ a^* & 1 & c^*
\end{array}}, $$
has eigenvalues $0$, $\pm \omega i$, where $\omega=\sqrt{\frac{1-b^2d^3}d}$.
For the critical values, system (\ref{FitzNag}) can be written as
$$\vetre{\dot x}{\dot y}{\dot z}=A\vetre xyz + F(x)\vetre 001,$$
where
$F(x) = x^3 - (a^*+1) x^2$, and  perform the linear transformation 
\beq \label{CLFitzNag}
\vetre xyz = P\vetre XYZ = \parent{ \begin{array}{ccc} 1 & 0 & d \\  b^2 d^2 &  b d \omega & 1 \\ 0 & -\omega & 0
\end{array}} \vetre XYZ,
\eeq
and the  time-reparametrization $dT=\omega dt$. Then, we obtain the system 
$$\vetre{X'}{Y'}{Z'}=\parent{ \begin{array}{cc|c} 0 & -1 &  \\  1 & 0 &   \\  \hline &   &  0
\end{array}} \vetre XYZ + F(X+dZ) \vetre{k_1}{k_2}{k_3},$$
where $\vetret{k_1}{k_2}{k_3} = P^{-1}\vetret 00{1/\omega}$, i.e.:
$$
k_1 = -b d / \omega^3, \
k_2 = -1/\omega^2, \ 
k_3 = b /  \omega^3.
$$

First, we reduce the principal part of the above system according to Proposition \ref{finCorl}, by means of the near-identity transformation 
\beq \label{NITFitzNag}
\hxx = X, \hyy = Y ,\hzz = Z +  \fracp{k_3 (d-1) }{ 2 \omega d}   X Y. 
\eeq
Then, we obtain
$$\vetrer{\hxx'}{\hyy'}{\hzz'}=
\vetrer{-\hyy}{\hxx}{\frac C2(\hxx^2+\hyy^2)}+
\textrm{HOT}
$$
where 
$C=k_3(1-d)/d$
and HOT stands for higher-order terms with respect to the type $\ttr$.
Next, we use the rescaling
\beq \label{RescalFitzNag}
\tx= \hxx, \ty= \hyy, \tz=\fracp 1{C}  \hzz.
\eeq
In this way, we obtain a system 
$$\vetrer{\tx'}{\ty'}{\tz'}=
\vetrer{-\ty}{\tx}{\frac 12 (\tx^2+\ty^2)}+
\textrm{HOT},
$$
whose principal part agrees with those given in (\ref{princpartHZ}).

The expression for the normal form coefficients, for the critical parameter values (see Theorem \ref{teoHZFN}), are:
$$
a_1 
=  -d^2/(2 \omega)
,\ 
b_1 
= d^2/(2\omega)
, \ 
c_1 
= -d / ( 2b )
.
$$

To compute a  parametric normal form  for the Hopf-zero singularity, as indicated in Theorem \ref{teoHZFN}, we 
consider system (\ref{FitzNag}) with parameters values 
$$a\approx a^*=-1/d, \ c\approx c^*=bd, \ d(1-b^2d^3)>0.$$
Let us consider the small parameters  $\Delta_a=a-a^*$, $\Delta_c=c-c^*$.
For parameter values close to the critical ones, system (\ref{FitzNag}) can be written as 
$$\vetre{\dot x}{\dot y}{\dot z}=A\vetre xyz +\Delta_A\vetre xyz + F(x)\vetre 001 + \Delta_F(x)\vetre 001 ,$$
where
$$\Delta_A=\parent{ \begin{array}{ccc} 0 & 0 & 0 \\  0 & 0 & 0 \\ \Delta_a & 0 & \Delta_c
\end{array}}, \textrm{ and } \ \Delta_F(x) = - \Delta_a  x^2.$$
Performing the linear transformation  (\ref{CLFitzNag})
and the  time-reparametrization $dT=\omega dt$, we obtain the system 
\bean
\vetre{X'}{Y'}{Z'}&=&\parent{ \begin{array}{cc|c} 0 & -1 &  \\  1 & 0 &   \\ \hline   &   &  0 \end{array}} \vetre XYZ 
+ F(X+dZ) \vetre{k_1}{k_2}{k_3}  \\ &&
+ \fracp 1{\omega^3} \parent{ \begin{array}{ccc}   -b d  \Delta_a &  b d  \omega \Delta_c  & -b d^2  \Delta_a \\  - \omega \Delta_a   &  \omega^2 \Delta_c   &  -  d \omega \Delta_a \\  b  \Delta_a &  -b \omega \Delta_c  &  b d  \Delta_a   \end{array}} \vetre XYZ 
+ \Delta_F(X+dZ) \vetre{k_1}{k_2}{k_3} .
\eean
Hence, according to the notation of the proof of Proposition \ref{ProposMiniversal}, we have
\bean
&&
\alpha_1 = - b d/\omega^3 \Delta_a,\
\alpha_2 = - b d/\omega^2 \Delta_c,
\\ &&
\beta_1 = 1/\omega \Delta_c,\
\beta_2 = -1/\omega^2 \Delta_a,
\\ &&
\gamma_1 = b/\omega^3 \Delta_a,\
\gamma_2 = - b/\omega^2 \Delta_c,\
\gamma_3 = b d/\omega^3 \Delta_a,\
\gamma_4 = 0,\
\gamma_5 = k_1 \Delta_a,\
\gamma_6 = 0.
\eean

 By means of the transformations indicated in the proof of Proposition \ref{ProposMiniversal}, we obtain the expressions for the unfolding parameters $\epsilon,\delta$, that are given by
 \bean 
 \epsilon &=& \frac{ d^2 (1-d^2)  (b d - 2 \omega^3)}{\omega^3 (1+d^2)^2 }
 \Delta_a
-
\frac{ (1-d^2)}{\omega (1+d^2)^2 }
\Delta_c
 + {\cal O}\left(\left|\Delta_a,\Delta_c\right|^2\right).
\\
\delta &=& \frac{bd(d^2-1)}{\omega^3(1+d^2)} \Delta_a  
+ {\cal O}\left(\left|\Delta_a,\Delta_c\right|^2\right).
\eean
Recall that $\Delta_a=a+\frac 1d$ and $\Delta_c=c-bd$ are small parameters.

The parametric normal form (\ref{NF_HoZe}) can be used to predict the presence of periodic behavior 
in the Hopf-zero bifurcation, as shown in the following theorem.

\begin{theorem}\label{teoOPHZFN} 
Let us consider the parametric normal form (\ref{NF_HoZe}), and assume that the normal form coefficient $a_1$ is non-zero.
Then, for $\epsilon$, $\delta$ small enough satisfying 
$\epsilon\parent{a_1\delta-b_1\epsilon}>0$, there exists a periodic orbit, which is given by
\bean
x(t)&=& \frac{\sqrt2}{|a_1|} \sqrt{\epsilon\parent{a_1\delta-b_1\epsilon}} \,\cos\parent{1-\frac{c_1\epsilon}{a_1}}+ {\cal O}\left(\left|\epsilon,\delta\right|^2\right),\\
y(t)&=& \frac{\sqrt2}{|a_1|} \sqrt{\epsilon\parent{a_1\delta-b_1\epsilon}} \,\sin\parent{1-\frac{c_1\epsilon}{a_1}}+ {\cal O}\left(\left|\epsilon,\delta\right|^2\right),\\
z(t)&=&- \frac{\epsilon}{a_1} \sqrt{\epsilon\parent{a_1\delta-b_1\epsilon}}+ {\cal O}\left(\left|\epsilon,\delta\right|^2\right).
\eean
Moreover, the characteristic multipliers of such periodic orbit are given, in first approximation, by
$$
\llave{\exp\parent{{i \frac{a_1\delta-2b_1\epsilon\pm\sqrt{(a_1\delta-2b_1\epsilon)^2+8a_1\epsilon(a_1\delta-b_1\epsilon)}}{2a_1} }},1}
.$$
\end{theorem}

\begin{proof}
It is enough to introduce cylindrical coordinates in system (\ref{NF_HoZe}), next compute the equilibrium point located outside the radial axis and finally 
take into account that this equilibrium corresponds to a periodic orbit in system (\ref{NF_HoZe}).
\end{proof}

In the case of the FitzHugh-Nagumo  system (\ref{FitzNag}), the above theorem can be applied because $b\neq 0$ 
(this implies $a_1\neq 0$). 
Moreover, the condition $\epsilon\parent{a_1\delta-b_1\epsilon}>0$ for the existence of the periodic orbit reads as
$$ (1/2)d^2(d^2-1)^2 
\parent{ (b+2d\omega^3) d\Delta_a  + \omega^2\Delta_c }
\parent{ (bd-2\omega^3)d^2\Delta_a  - \omega^2\Delta_c }
/
(\omega^7(1+d^2)^4)>0.$$
Recalling that $\omega=\sqrt{\frac{1-b^2d^3}d}$, we obtain that this condition is fulfilled if, and only if,  
$$
d\neq \pm 1 \ 
\textrm{ and }
\
\textrm{sign}\parent{ (b+2d\omega^3) d\Delta_a  + \omega^2\Delta_c }=\textrm{sign}
\parent{ (bd-2\omega^3)d^2\Delta_a  - \omega^2\Delta_c }.$$ 

\

\noindent\textbf{Acknowledgements.}
{This research was partly supported by the
\emph{Ministerio de  Ciencia e Innovaci\'on, fondos FEDER}
(projects MTM2014-56272-C2-X-P) and by the \emph{{Consejer\'{\i}a de Econom\'{\i}a, Innovaci\'on, Ciencia y
Empleo} de la Junta de Andaluc\'{\i}a} (projects
{P12-FQM-1658}, TIC-130, FQM-276).
}

\end{document}